\documentclass[10pt]{article}
\usepackage[a4paper, top=1.5cm, bottom=1.5cm, left=1.5cm, right=1.5cm]{geometry}

\usepackage{amsmath,amssymb}

\usepackage{changepage}

\usepackage{textcomp,marvosym}

\usepackage[hidelinks]{hyperref}
\usepackage{nameref}

\usepackage[nopatch=eqnum]{microtype}
\DisableLigatures[f]{encoding = *, family = * }

\usepackage[table]{xcolor}

\usepackage{array}
\usepackage{booktabs}

\usepackage[aboveskip=1pt,labelfont=bf,labelsep=period,justification=raggedright,singlelinecheck=off]{caption}

\makeatletter
\renewcommand{\@biblabel}[1]{\quad#1.}
\makeatother

\usepackage{subcaption}
\usepackage{tikz}
\usepackage{xcolor}
\usepackage{colortbl}
\usepackage{graphics}
\usepackage{graphicx}
\usepackage{pgfplots}
\usepackage{parboxx}
\usepackage{siunitx}
\usepackage{float}
\usepgfplotslibrary{groupplots}
\usepackage{pgfplotstable}
\usepgfplotslibrary{statistics}
\usepackage[acronym]{glossaries}
\pgfplotsset{compat=1.18}
\usetikzlibrary{arrows.meta, positioning, fit, shapes.multipart}

\usepackage{lastpage,fancyhdr,graphicx}
\usepackage{epstopdf}
\fancyheadoffset[L]{2.25in}
\fancyfootoffset[L]{2.25in}
\newcommand{\myVec}[1]{\mathbf{#1}}
\newcommand{\myTens}[1]{\mathbf{#1}}
\newcommand{\alphaEpiRV}{\alpha_{\mathrm{Epi,RV}}}
\newcommand{\alphaEndoRV}{\alpha_{\mathrm{Endo,RV}}}
\newcommand{\alphaEndoLV}{\alpha_{\mathrm{Endo,LV}}}
\newcommand{\alphaEpiLV}{\alpha_{\mathrm{Epi,LV}}}
\newcommand{\alphaRBM}{\alpha_{\mathrm{RBM}}}

\newcommand{\gammaEpiRV}{\gamma_{\mathrm{Epi,RV}}}
\newcommand{\gammaEndoRV}{\gamma_{\mathrm{Endo,RV}}}
\newcommand{\gammaEndoLV}{\gamma_{\mathrm{Endo,LV}}}
\newcommand{\gammaEpiLV}{\gamma_{\mathrm{Epi,LV}}}
\newcommand{\gammaRBM}{\gamma_{\mathrm{RBM}}}

\newcommand{\refL}{\myVec{e}_{\mathrm{l}}}
\newcommand{\refT}{\myVec{e}_{\mathrm{t}}}
\newcommand{\refN}{\myVec{e}_{\mathrm{n}}}
\newcommand{\fLN}{\myVec{f}_{\mathrm{LN}}}

\newcommand{\refTPrime}{\myVec{e}'_{\mathrm{t}}}
\newcommand{\refNPrime}{\myVec{e}'_{\mathrm{n}}}

\newcommand{\boardApex}{\Gamma_{\mathrm{Apex}}}
\newcommand{\boardBase}{\Gamma_{\mathrm{Base}}}
\newcommand{\boardEndoRV}{\Gamma_{\mathrm{Endo,RV}}}
\newcommand{\boardEndoLV}{\Gamma_{\mathrm{Endo,LV}}}
\newcommand{\boardEpi}{\Gamma_{\mathrm{Epi}}}

\newcommand{\ActTension}{{T}_{\mathrm{a}}}
\newcommand{\mF}{m_{\mathrm{f}}}
\newcommand{\mS}{m_{\mathrm{s}}}
\newcommand{\mN}{m_{\mathrm{n}}}
\newcommand{\sigmaF}{\sigma_{\mathrm{f}}}
\newcommand{\sigmaS}{\sigma_{\mathrm{s}}}
\newcommand{\sigmaN}{\sigma_{\mathrm{n}}}
\newcommand{\CBstiffness}{a_{\mathrm{XB}}}

\newcommand{\SF}{\mathrm{SF}}
\newcommand{\SL}{\mathrm{SL}}
\newcommand{\Calcium}{[\mathrm{Ca}^{2+}]}

\makeglossaries
\newacronym{moi}{MOI}{Mesoscopic Optical Imaging}
\newacronym{rbm}{RBM}{Rule-Based Method}
\newacronym{rbms}{RBMs}{Rule-Based Methods}
\newacronym{ldrbm}{LDRBM}{Laplace-Dirichlet Rule Based Method}
\newacronym{ldrbms}{LDRBMs}{Laplace-Dirichlet Rule Based Methods}
\newacronym{frs}{FRS}{Fiber Reference System}
\newacronym{mrf}{MRS}{Myocadium Reference System}
\newacronym{lv}{LV}{Left Ventricle}
\newacronym{rv}{RV}{Right Ventricle}
\newacronym{la}{LA}{Left Atrium}
\newacronym{ra}{RA}{Right Atrium}
\newacronym{edp}{EDP}{End-Diastolic Pressure}
\newacronym{edv}{EDV}{End-Diastolic Volume}
\newacronym{esp}{ESP}{End-Systolic Pressure}
\newacronym{esv}{ESV}{End-Systolic Volume}
\newacronym{fe}{FE}{Finite Elements}
\newacronym{ef}{EF}{Ejection Fraction}
\newacronym{sf}{$\SF$}{Stress Factor}
\newacronym{qoi}{QoI}{Quantity of Interest}
\newacronym{dtmri}{DTMRI}{Diffusion Tensor Magnetic Resonance Imaging}
\newacronym{bc}{BCs}{Boundary Conditions}

\DeclareSIUnit\mmHg{mmHg}

\begin{document}
\vspace*{0.2in}

\begin{flushleft}
{\Large
\textbf{The functional impact of myofiber macroscopic organization and disarray in computational models of the murine heart} 
}
\newline
\\
Carlo Guastamacchia\textsuperscript{1},
Roberto Piersanti\textsuperscript{2},
Francesco Giardini\textsuperscript{3},
Raffaele Coppini\textsuperscript{4},
Cecilia Ferrantini\textsuperscript{4},
Luca Dede'\textsuperscript{1},
Leonardo Sacconi\textsuperscript{3,5},
Francesco Regazzoni\textsuperscript{1,*},
\\
\bigskip
\textbf{1} MOX - Department of Mathematics, Politecnico di Milano, Milan, Italy.
\\
\textbf{2} Department of Theoretical and Applied Sciences
(DiSTA), Università degli Studi eCampus, Novedrate, Italy.
\\
\textbf{3} Institute for Experimental Cardiovascular
Medicine, University Heart Center Freiburg – Bad
Krozingen, Medical Faculty and Medical Center –
University of Freiburg, Freiburg im Breisgau,
Germany.
\\
\textbf{4} Department of Experimental and Clinical
Medicine, University of Florence, Florence, Italy.
\\
\textbf{5} Institute of Clinical Physiology, National Research
Council (IFC-CNR), Florence, Italy.
\bigskip

%
%


\textcurrency Current Address: Piazza Leonardo da Vinci 32, 20133 Milano, Italia 

* francesco.regazzoni@polimi.it

\end{flushleft}


%
\section*{Abstract}
A major challenge in computational models of cardiac electromechanics is the reconstruction of myocardial fiber architecture, as direct in vivo measurements of fiber orientation are not feasible. Consequently, rule-based methods are commonly adopted as surrogates, relying on empirical descriptions of fiber organization combined with patient-specific geometries.
This study investigates the respective roles of macroscopic fiber architecture and microscopic fiber disarray in cardiac electromechanical simulations. A high-fidelity biventricular electromechanical model of a murine heart was developed using a high-resolution myocardial fiber field obtained via mesoscopic optical imaging, which serves as a reference ground truth. A spatial smoothing strategy is introduced to decouple macroscopic fiber organization from local disarray, and the resulting responses are also compared with those obtained using a rule-based fiber field.
The results show that passive mechanics and electrophysiological activation are only weakly affected by fiber disarray, with global chamber compliance and activation times remaining largely unchanged across different fiber descriptions. In contrast, active mechanics is highly sensitive to fiber architecture. Moderate regularization of the experimentally measured fiber field enhances the ventricular pumping efficiency of the computational model by reducing microscopic disarray while preserving the macroscopic helical organization, whereas excessive smoothing or rule-based fiber reconstructions lead to unphysiologically strong or inefficient contraction.
Within this framework, two commonly adopted surrogate strategies to account for fiber disarray are investigated: (i) a reduction of the effective cross-bridge stiffness in the active tension model, and (ii) the introduction of controlled misalignment between active tension and the local fiber direction. While both approaches reproduce global hemodynamic indicators comparable to the reference case, an effective reduction of contractility -- despite its phenomenological nature -- provides a closer match to the reference strain patterns than the introduction of orthogonal active stress components. Overall, the results highlight the dominant role of macroscopic fiber architecture in active mechanics and reveal important limitations of commonly adopted surrogate approaches for modeling fiber disarray.

\section*{Author summary}
To integrate cardiac computational models into routine clinical practice, it is essential to validate the modeling strategies proposed in the scientific literature. A major challenge in this context is the reconstruction of myocardial fiber fields, as obtaining precise in-vivo measurements of fiber orientation within patient-specific organs remains infeasible. An additional source of complexity arises from fiber disarray, which refers to local deviations of fiber orientation from the average direction.

Current state-of-the-art methods for reconstructing cardiac fiber fields from heart geometry typically neglect the presence of fiber disarray. In this work, we analyze a dataset containing high-resolution fiber information obtained from ex-vivo measurements of a mouse heart. We introduce a methodology to assess the relative impact of the mean fiber direction and fiber disarray on electromechanical simulations, and we evaluate the reliability of a state-of-the-art fiber surrogation approach.

Our results indicate that fiber orientation has only a limited effect on electrical signal propagation and on the passive mechanical properties of the heart. However, when active tension generated by myocardial fibers is included in the model, fiber orientation and disarray have a strong influence on cardiac contraction and local tissue deformations. urthermore, the results show that commonly used surrogate fiber models, even when combined with specific modeling techniques to account for disarray, are not
sufficient to accurately reproduce the local effects associated with the real fiber field.

\newcommand{\TableBaseline}[2]{%
\begin{table}[!htbp]
\centering
\scriptsize
\begin{tabular}{lllll}
	\toprule
	\textbf{Component} & \textbf{Model Parameter} & \textbf{Value} & \textbf{Method parameter} & \textbf{Value} \\
\midrule
	\textbf{Electrophysiology} & Electrophysiology approach & Eikonal & FE space degree & 1 \\
 & Longitudinal conductivity $\sigmaF$ & \SI{2.0e-4}{\meter^2 \second^{-1}} & Algorithm & Continuation \\
 & Depolarization coefficient & \SI{55}{\second^{-1/2}} & Initial ramp step & 0.1 \\
 & Transversal conductivity $\sigmaS$ & \SI{1.0e-4}{\meter^2 \second^{-1}} & Step factor after failure & 0.5 \\
 & Cross-fiber conductivity $\sigmaN$& \SI{1.0e-4}{\meter^2 \second^{-1}} & Step factor after success & 1.2 \\
 & &  & Maximum ramp increment & 1 \\
 &  &  & Minimum ramp increment & 0 \\
 & & & Linear solver & GMRES, \\
 & & & & Maximum iterations =1000, \\
 & & & & Tolerance=1e-15,\\
\addlinespace
	\textbf{Fibers} & Regularization radius $\ell$ & \SI{0.0}{\milli \meter} & Element type & Tet \\
\addlinespace
	\textbf{Mechanics} & Density & \SI{1e3}{\kilogram \meter^{-3}} & FE space degree & 1 \\
 & Viscosity & \SI{0}{\pascal\second} & Pericardium Robin BC & $K_{\text{norm}}=\SI{30e5}{\meter^{-1} \pascal }$, \\
 & Fiber active stress factor $\mF$& 1.0 & & $K_{\text{tan}}=\SI{30e4}{\meter^{-1} \pascal }$, \\
 & Sheet active stress factor $\mS$ & 0.0 & & $C_{\text{norm}}=\SI{6e4}{\second\meter^{-1}\pascal}$, \\
 & Cross active stress factor $\mN$ & 0.0 & & $C_{\text{tan}}=\SI{6e3}{\second\meter^{-1} \pascal}$. \\
 & & & Base Robin BC & $K_{\text{norm}}=\SI{2e5}{\meter^{-1} \pascal}$, \\
 & Usyk parameters & $b_{\text{f}}=5.0$, $b_{\text{s}}=3.7$, & & $K_{tan}=\SI{2e5}{\meter^{-1} \pascal}$, \\
 & & $b_{\text{n}}=3.7$, $b_{\text{fs}}=2.6$, & & $C_{\text{norm}}=\SI{2e4}{\second \meter^{-1} \pascal}$, \\
 & & $b_{\text{fn}}=3.7$, $b_{\text{sn}}=3.7$, & & $C_{\text{tan}}=\SI{2e3}{\second \meter^{-1} \pascal}$. \\
 & & $\text{bulk}=\SI{5e4}{\pascal}$, $c=\SI{2e3}{\pascal}$ & Non-linear solver & Maximum iterations =50, \\
 & & & & Absolute tolerance=1e-9, \\
 & & & & Relative tolerance=1e-9 \\
 & & & Linear solver & GMRES, \\
 & & & & Maximum iterations=1000, \\
 & & & & Tolerance=1e-8, \\
\addlinespace
	\textbf{Activation} & LV contractility ratio & 1.0 & FE space degree & 1 \\
 & RV contractility ratio & 0.5 & Time step (0D) & \SI{2e-5}{\second} \\
 & Sarcomere slack length & $\SI{1.9}{\micro \meter}$ & Time step RU update & \SI{2.0e-5}{\second} \\
 & Model & Activation-RDQ20-MF & Linear solver & GMRES, \\
 & Sarcomere initial length & \SI{2.2}{\micro \meter} & &Maximum iterations=1000, \\
 & Calcium initial concentration & \SI{0.1656}{\micro \mole \liter^{-1}} & & Tolerance=1e-10,\\
 & Initialization time & \SI{0.2}{\second} & & \\
 & Sarcomere parameters & $LA=\SI{1.25}{\micro \meter}$, $LM=\SI{1.65}{\micro \meter}$, & & \\
 & & $LB=\SI{0.18}{\micro \meter}$,  $\mu_{0, \text{fP}}=\SI{32.68}{\second^{-1}}$, & & \\
 & & $SL_{0}=\SI{2.2}{\micro \meter}$, $Q=2$, & & \\
 & & \multicolumn{2}{l}{$K_{\mathrm{d},0}=\SI{0.65}{\micro \meter}$, $\alpha_{\text{Kd}}=\SI{-0.83}{\micro \mole \liter^{-1} \micro \meter^{-1}}$,} & \\
 & & $\mu=10$, $\gamma=12$, & & \\
 & & $K_{\mathrm{off}}=\SI{180}{\second^{-1}}$, $K_{\mathrm{basic}}=\SI{180}{\second ^{-1}}$, & & \\
 & & $r_0=\SI{134.3}{\second^{-1}}$, $\alpha=25.18$, & & \\
 & & $\mu_{1, \text{fP}} = \SI{0.77}{\second^{-1}}$, $\CBstiffness=\SI{15e6}{\pascal}$, & & \\
 & & $TnC_{\mathrm{max}}=\SI{70}{\micro \mole \liter^{-1}}$ & & \\
\addlinespace
	\textbf{Circulation} & Heartbeat period & \SI{0.2}{\second} & Numerical scheme & Explicit Euler \\
 & Right atrium & \multicolumn{2}{l}{Elastance $A=\SI{1000}{\mmHg \milli \liter^{-1}}$,} & \\
 & & $B=\SI{40}{\mmHg \milli \liter^{-1}}$, $V_0=\SI{0.0008}{\milli \liter}$, & & \\
 & & $\text{Duration of contraction} = 0.38$, & & \\
 & & $\text{Duration of relaxation} = 0.38$, & & \\
 & & $t_0=0.83$, & & \\
 & Pulmonary arterial system & \multicolumn{2}{l}{$R=\SI{5.216}{\mmHg \second \milli\liter^{-1}}$, $C=\SI{0.0024}{\milli \liter \mmHg^{-1}}$,} & \\
 & & \multicolumn{2}{l}{$R_{\mathrm{up}}=\SI{0.0}{\mmHg \second \milli\liter^{-1}}$, $L=\SI{2064e-4}{\mmHg \second^2 \milli \liter^{-1}}$}& \\
 & Pulmonary venous system & \multicolumn{2}{l}{$R=\SI{55.45}{\mmHg \second \milli\liter^{-1}}$, $C=\SI{0.0387}{\milli \liter \mmHg^{-1}}$} & \\
 & & $L=\SI{2064e-4}{\mmHg \second^2 \milli \liter^{-1}}$ & & \\
 & Left atrium & \multicolumn{2}{l}{Elastance $A=\SI{140}{\mmHg \milli \liter^{-1}}$, $B=\SI{360}{\mmHg \milli \liter^{-1}}$,} & \\
 & & $\text{Duration of contraction}=0.36$, & \\
 & & $\text{Duration of relaxation}=0.36$, & \\
 & & \multicolumn{2}{l}{$t_0=0.85$, $V_0=\SI{0.0008}{\milli \liter}$} & \\
 & Systemic arterial system & \multicolumn{2}{l}{$R=\SI{500}{\mmHg \second \milli\liter^{-1}}$, $C=\SI{0.0008}{\milli \liter \mmHg^{-1}}$,} & \\
 & & \multicolumn{2}{l}{$R_{\mathrm{up}}=\SI{0.0}{\mmHg \second \milli\liter^{-1}}$, $L=\SI{2064e-3}{\mmHg \second^2 \milli \liter^{-1}}$,} & \\
 & Systemic venous system & \multicolumn{2}{l}{$R=\SI{400}{\mmHg \second \milli\liter^{-1}}$, $C=\SI{0.0145}{\milli \liter \mmHg^{-1}}$,} & \\
 & & $L=\SI{2064e-4}{\mmHg \second^2 \milli \liter^{-1}}$ & & \\
 & Valves & \multicolumn{2}{l}{$R_{\mathrm{min}}=\SI{7.5}{\mmHg \second \milli\liter^{-1}}$, $R_{\mathrm{max}}$=\SI{75000}{\mmHg \second \milli\liter^{-1}},} & \\
\bottomrule
\end{tabular}
\caption{#2}
\label{tab:#1}
\end{table}
}

\newcommand{\TuningMap}[2]{
\begin{figure}[t]
\centering
\begin{tikzpicture}[
    font=\sffamily,
    node distance=1.4cm and 2.2cm,
    box/.style={draw, rounded corners, thick, align=center, text width=3cm, minimum height=1.5cm, fill=white},
    title/.style={font=\bfseries\large, align=center},
    arrow/.style={-{Stealth[length=3mm]}, thick},
]

\node[title] (col1) {Blood Circulation};
\node[title, right=1.6cm of col1] (col2) {Mechanics};
\node[title, right=1.6cm of col2] (col3) {Electrophysiology};
\node[title, right=1.6cm of col3] (col4) {Activation};

\node[box, below=0.5cm of col1] (bc1) {Full 0D circulation\\tuning on literature QoI};
\node[box, below=1.0cm of bc1] (bc2) {0D emulator};
\node[box, below=1.0cm of bc2] (bc3) {Parametric 0D emulator};

\node[box, right=1.2cm of bc1] (mech1) {Usyk’s parameters tuning to obtain PV curve with $Ta = \SI{0}{\mega \pascal}$};
\node[box, right=1.2cm of bc2] (mech2) {Full electromechanical simulation};
\node[box, right=1.2cm of bc3] (mech3) {$\CBstiffness$ tuning};

\node[box, right=1.2cm of mech1] (elec1) {Calibration to obtain depolarization speed from experimental data};

\node[box, right=1.2cm of elec1] (act1) {Calibration of the cell model};

\draw[arrow, orange!80] (mech1.south) -- node[midway, right]{3} (mech2.north);
\draw[arrow, yellow!80] (elec1.south west) -- node[midway, below]{1} (mech2.north east);
\draw[arrow, red!80] (act1.south) |- node[midway, below]{4} (mech2.east);
\draw[arrow, green!80] (mech2.west) -- node[midway, above]{5} (bc2.east);
\draw[arrow, green!80] (bc2.north) -- node[midway, right]{6} (bc1.south);
\draw[arrow, green!80] (bc1.south east) -- node[midway, left]{\rotatebox{-43}{2 and 7}} (mech2.north west);
\draw[arrow, blue!80] (mech2.south west) -- node[midway, above]{8} (bc3.north east);
\draw[arrow, blue!80] (bc3.east) -- node[midway, above]{9} (mech3.west);
\draw[arrow, blue!80] (mech3.north) -- node[midway, right]{10} (mech2.south);
\end{tikzpicture}
\caption{#2}
\label{fig:#1}
\end{figure}
}
\section{Introduction}\label{sec:Introduction}
Cardiac computational models are increasingly being used in modern medical practice \cite{peirlinck2021precision, niederer2019computational}.
These models have demonstrated the ability to reproduce action potential propagation under both physiological \cite{trayanova2024computational, vazquez2011massively} and pathological conditions \cite{pagani2021computational, stella2020integration, mehri2025multi, wulfers2024whole}, to simulate mechanical contraction in electromechanical frameworks \cite{trayanova2011electromechanical, trayanova2011cardiac, griffith2013electrophysiology}, and to model blood flow in fluid dynamic simulations \cite{fumagalli2022image}.
By leveraging these models, it is possible to investigate and predict specific quantities of interest \cite{montino2024personalized}, prognosticate disease evolution \cite{gao2017estimating}, and test new medical treatments in silico \cite{rodero2023systematic}.
The ultimate goal in this field is the development of patient-specific digital twins capable of delivering personalized analyses and therapeutic strategies \cite{fumagalli2024role, trayanova2012computational, viola2023gpu, sack2018construction, taylor2009patient, chapelle2009numerical, corti2022impact}.
Achieving this objective requires an accurate representation of the heart’s muscles structure and function; however, with current technological capabilities, such detailed reconstruction remains unfeasible on a patient specific basis \cite{piersanti2025defining}.

The cardiac muscle is composed of cardiomyocytes arranged into fibers, commonly referred to as myofibers. As described by \cite{legrice1995laminar}, the myocardial fiber architecture consists of an orderly laminar organization of myofibers, characterized by extensive cleavage planes separating adjacent muscle layers. 
In transmural sections, these planes extended radially from the endocardium (the inner surface) to the epicardium (the outer surface), aligning with the local myofiber orientation observed in tangential sections. This well-organized laminar architecture can be described in terms of three material symmetry axes: the logitudinal, the sheet, and the normal directions. The longitudinal axes of myocytes, which constitute laminae, show a well-defined helical organization that progressively rotates along the transmural axis from the endocardium to the epicardium \cite{streeter1969fiber}.
However, the complexity of the myofiber architecture increases when considering the heterogeneity in macroscopic fiber orientation and the presence of physiological disarray \cite{lombaert2012human}. The disarray of fibers consists in a local misalignment of the fibers with respect to the mean direction \cite{teare1958asymmetrical}.

The characteristic myofiber configuration has a major impact on the heart function. Fiber orientation affects action potential propagation within the muscle, as the spread of epicardial excitation is considerably 
faster parallel to the longitudinal axes of the cardiac fibers than perpendicular to them \cite{roberts1979influence}. 
Furthermore, fiber disarray can alter tissue electrical conductivity, leading to a more isotropic propagation of the wavefront \cite{punske2005effect}. 
In addition, the tissue’s mechanical contraction, driven by the propagation of electrical signals, is highly dependent on the alignment of the muscle fibers \cite{guccioneII}. Due to its anisotropic nature, cardiac muscle exhibits direction-dependent material stiffness determined by the local myofiber architecture along the three principal \cite{guccione1991finite}. Morover, cardiomyocytes are also responsible for active contraction, producing force mainly aligned with the fiber orientation \cite{guyton2006text}. Since myofibers play a pivotal role in cardiac function -- affecting electrophysiological behavior, passive mechanical properties, and active contraction -- precise representation of their spatial architecture is crucial in cardiac computational modeling. In particular, cardiac digital twins require to represent patient-specific fiber architectures \cite{gerach2021electro, piersanti2025defining}. Furthermore, assessing the impact of fiber and disarray is key to unraveling the complex pathophysiological mechanisms underlying the cardiac diseases, such as hypercontractility, hypocontractility, and fibrosis \cite{mehri2025multi}.

The current \textit{de facto} standard imaging technique for reconstructing the fiber architecture is \acrfull{dtmri} \cite{hsu1998magnetic,helm2005ex}, 
which has achieved resolutions of \SI{400}{\micro\meter} in ex-vivo human hearts \cite{pashakhanloo2016myofiber} and \SI{43}{\micro\meter} in ex-vivo mouse hearts \cite{angeli2014high}. 
However, DTMRI suffers from poor signal-to-noise ratio and long acquisition times \cite{teh2016resolving}. As an alternative, micro-computed tomography \cite{gonzalez2017whole}, which does not require long preparation or measurement times, has reached a resolution of \SI{10}{\micro \meter} in ex-vivo rat hearts. More recently, the authors of \cite{walsh2021imaging} demonstrated that, by means of Hierarchical Phase-Contrast Tomography (HiP-CT), it is possible to reconstruct an entire human heart at an isotropic resolution of approximately $\SI{25}{\micro\meter}$ per voxel, while integrating hierarchical local scans down to $\sim\SI{2}{\micro\meter}$ per voxel for the detailed analysis of cardiac musculature and cardiomyocyte organization. Nevertheless, the applicability of such approaches is limited by the high cost of synchrotrons for X-ray production \cite{wang2024investigating} which severely restricts their accessibility.
Another imaging technique is shear-wave imaging, which achieved a resolution of \SI{200}{\micro\meter} in ex-vivo porcine hearts \cite{lee2011mapping}. 
It is worth noting that optical imaging techniques can achieve micrometric or even sub-micrometric resolution; however, this is typically limited to thin tissue sections or optically cleared samples with constrained thickness \cite{dileep2023cardiomyocyte}, and does not allow imaging of the intact whole organ at comparable resolution.
To overcome the aforementioned limitations, recent advances in optical imaging techniques \cite{tolstik2024cardiac} have enabled mesoscale reconstruction of cardiac anatomy at the whole-organ level. In particular, the combination of tissue clearing methods with a new generation of mesoSPIM microscopes enables whole-heart reconstruction at an isotropic resolution of \SI{3.25}{\micro \meter} × \SI{3.00}{\micro \meter} in ex vivo mouse hearts, as demonstrated by  \cite{giardini2025correlative}.  All of the imaging techniques mentioned above are performed ex vivo. Currently, in vivo fiber identification remains limited by their relatively coarse spatial resolution \cite{toussaint2013vivo, froeling2014diffusion,nguyen2016vivo,nielles2014erratum}. Therefore, contemporary fiber-imaging techniques are largely impractical for building patient-specific cardiac computational models.

Given the challenges associated with obtaining patient-specific fiber fields, mathematical models, known as \acrfull{rbms}  are commonly employed in cardiac computational models. \acrshort{rbms} surrogate the characteristic myocardial fiber architecture. \acrshort{rbms} represent fiber orientations exploiting mathematically sound rules informed by histological or imaging data together with patient-specific cardiac geometry \cite{piersanti2021modeling, piersanti2025defining, bayer2012novel, doste2019rule, rossi2014thermodynamically, wong2014generating}.
The current state-of-the-art of \acrshort{rbms} is represented by the \acrfull{ldrbms}, which determine the myofiber direction by solving suitable Laplace-Dirichlet problems. Despite their widespread use, the impact of these methods on the reliability of models has yet to be definitively established. In \cite{doste2019rule}, the authors compared the fiber fields generated by ventricular \acrshort{ldrbm} with those obtained from diffusion tensor magnetic resonance imaging (\acrshort{dtmri}) at a spatial resolution of $\SI{400}{\micro \meter}$. The study reported a non-negligible angular discrepancy of approximately $\SI{30}{\degree}$. From an electrophysiological perspective, the works in \cite{bayer2012novel, doste2019rule} validated ventricular \acrshort{ldrbms} by comparing simulated activation maps with experimental measurements. The method proved effective in reproducing activation patterns across biventricular geometries. Specific \acrshort{rbms} have been developed for atria as they present more complex fiber architecture characterized by bundles with different mean fiber directions \cite{piersanti2021modeling, ferrer2015detailed, tobon2013three}. Nevertheless, the analysis of these models is out of the scope of this work. 

Still, should an exact representation of every single cardiomyocyte orientation be available, the computational cost of running a simulation of the full organ by resolving the cell scale would be unaffordable.
As a matter of fact, attempts to simulate the cardiac function at the cell scale remain limited to specimens composed by a few cells {\cite{huynh:BDDC,rosilhoBEM,goebel2025bddc}}.
Most of cardiac computational models are in fact homogenized, as they describe the tissue as a continuum, without resolving the single cells. 
In such models, fiber disarray is not explicitly represented, but it is implicitly incorporated either through effective parameters or through purposely defined corrective terms.
In particular, in electrophysiology models, fiber disarray is implicitly accounted for in the definition of the anisotropic diffusion tensor.
Similarly, in passive mechanics models, fiber disarray is encoded in the coefficients of the hyperelastic constitutive law \cite{guccione1991passive,holzapfel2009constitutive,nordbo2014computational}. Indeed, experimental tissue samples used for mechanical testing and models calibration inherently exhibit fibers with some degree of disarray \cite{ramadan2017standardized}. Conversely, in active mechanics, a dedicated modeling strategy is required  to account for the influence of fiber disarray \cite{niederer2006quantitative, land2017model, rice2008approximate}. To this end, the authors of \cite{usyk2000effect} proposed introducing cross-fiber activation to mimic the effects of fiber disarray within \acrshort{rbm} frameworks. 
In \cite{rossi2014thermodynamically}, \acrshort{rbms} were employed in biventricular electromechanical simulations that successfully reproduced realistic contraction patterns. In \cite{piersanti20223d}, biventricular electromechanical simulations incorporating cross-fiber active tension successfully reproduced experimental data of key mechanical biomarkers reported in the literature under physiological conditions. However, in both studies, no experimental ground truth was available for validation.
The authors in \cite{guan2020effect, guan2021modelling} developed a mechanical model of a biventricular geometry, comparing simulations based on \acrshort{rbm}-generated fibers with those using experimentally measured ones. The study showed that incorporating cross-fiber activation improved the agreement of pressure–volume (PV) loops with experimental data but failed to fully capture local deformation effects. However, this analysis, conducted on porcine heart model embedded with canine fibers, is limited by the absence of coupling with an electrophysiology model. In conclusion, to the best of our knowledge, no validation of \acrshort{rbms} using comprehensive electromechanical models and high-resolution fiber data has yet been reported in the literature.

Despite the non conclusive literature results, cross-fiber activation remains a widely adopted approach in standard electromechanical simulations employing \acrshort{rbms} \cite{fedele2023comprehensive, quarteroni2017integrated, usyk2002computational, piersanti20223d}. Moreover, while the sensitivity of mechanical \cite{palit2015computational} and electromechanical \cite{gil2019influence,eriksson2013influence,guan2020effect,piersanti20223d} simulations to fiber orientation has been shown in literature, the specific impact of fiber disarray with respect to the mean fiber orientation has, to the best of our knowledge, not yet been explored.

Motivated by these open issues, in this work we present a computational study based on a measured myofiber field obtained using mesoscopic optical imaging \cite{giardini2022mesoscopic} in a murine heart, enabling an electromechanical analysis at unprecedented spatial resolution. The different components of the model are calibrated using functional measurements acquired from the same specimen, including activation maps, while complementary information is drawn from the literature when direct measurements are not available.
To independently assess the functional impact of macroscopic fiber architecture and microscopic fiber disarray, we introduce a methodology based on spatial frequency decomposition that allows these two contributions to be disentangled. Leveraging the very high resolution and signal-to-noise ratio of the measured fiber field, together with the proposed disentangling approach, this study enables multiple analyses. First, we investigate the respective effects of macroscopic architecture and microscopic disarray on cardiac function, spanning passive mechanics, electrophysiology, and active contraction. Second, we provide a validation of commonly adopted modeling strategies used in practice to account for these features, namely \acrshort{rbm} on the one hand, and the inclusion of cross-fiber activation or modulation of cross-bridge stiffness on the other.

The remainder of the paper is organized as follows. 
In Sec.~\ref{sec:Models}, we present the methods used to characterize the macroscopic fiber architecture and local disarray, as well as the \acrshort{ldrbm} considered in this work.
In Sec.~\ref{sec:Baseline}, we describe the electromechanical model used to study the effect of fiber architecture on cardiac electromechanics.
In Sec.~\ref{sec:Results}, we present the results, while Sec.~\ref{sec:Discussion} is devoted to their discussion.
\section{Myofiber architecture: data and models}\label{sec:Models}

In this section, we introduce the notation employed to characterize the reference system and the angles defining the fibers architecture. Then we describe the experimental data and the procedure used to isolate the macroscopic fiber architecture from the microscopic disarray. Finally, we briefly describe the \acrshort{ldrbm} considered in this work.

\subsection{Fiber architecture and rotation angles}\label{subsec:fib_arch_angles}
\begin{figure}[t]
\includegraphics[width=\linewidth]{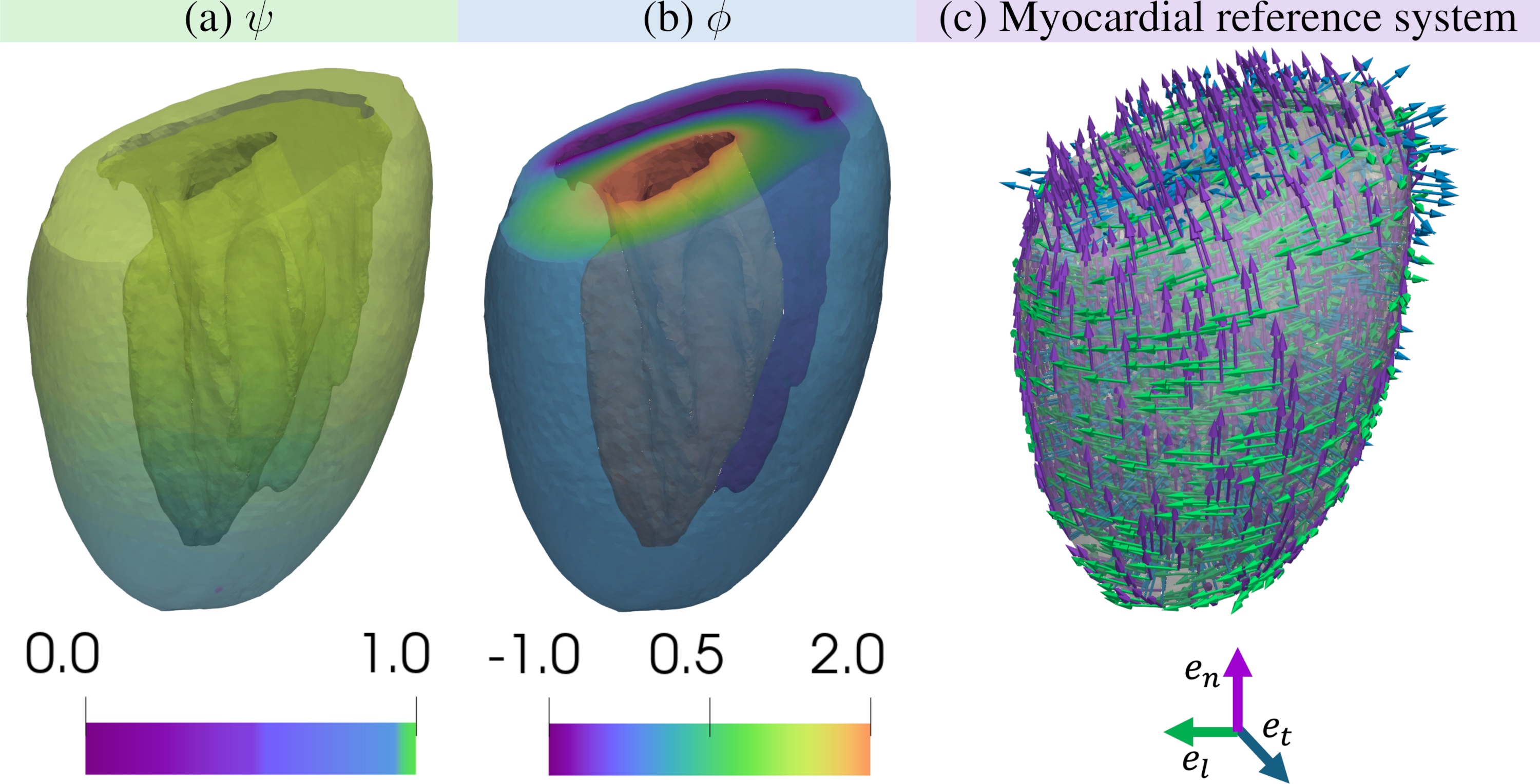}
\caption{\textbf{From myocardial to fiber reference .} (a) Rotation of $\refL$ by an angle $\alpha$ around the axis $\myVec{e}_t$ to obtain $\fLN$, followed by a rotation of $\fLN$ by an angle $\gamma$ in the plane $\text{span}\{\fLN, \myVec{e}_t \}$ (depicted in blue) to obtain $\myVec{f}$. (b) Rotation of $\myVec{e}_t$ by an angle $\gamma$ to obtain $\myVec{e}_t'$. (c) Rotation of  $\myVec{e}_t'$ by an angle $\beta$ around $\myVec{f}$ (in the lilac plane) to obtain $\myVec{s}$, and definition of $\myVec{n}$ as the normal to both $\myVec{f}$ and $\myVec{s}$.}
\label{fig:fig1}
\end{figure}

The fiber field is defined by three principal orthogonal directions, which characterize the tissue properties: $\myVec{f}$ is the longitudinal fiber direction; $\myVec{s}$ is the transmural sheet direction defining the orthogonal direction to the fiber and lying on the fiber's sheet plane; $\myVec{n}$ is the normal direction to the sheet plane. These three orthogonal directions compose the fiber reference system. 
In modeling practice, an additional triplet known as the myocardial reference system is often introduced \cite{piersanti2021modeling}. Its formulation relies solely on information about myocardial geometry and is designed to parametrize the transmural and apico-basal directions throughout the entire myocardium, thereby defining an orthotropic reference frame. Consequently, the myocardial reference frame can be used to analyze fiber orientation with respect to myocardial geometry or to surrogate the fiber field through \acrshort{rbms} \cite{piersanti2021modeling, piersanti2025defining}. The three orthogonal directions of the myocardial reference system are: the apico-basal direction $\refN$, going from the ventricular apex to the base; the transmural direction $\refT$, going from the endocardium to the epicardium; the circumferential direction $\refL$, defined as normal to the other two.
The fiber's principal directions ($\myVec{f}$, $\myVec{s}$, $\myVec{n}$) and the myocardial reference system ($\refL$, $\refT$, $\refN$) are linked by rotation angles. These allow to pass from one reference system to the other. In particular, three rotation angles are required to map the myocardial reference system into the fiber reference system \cite{agger2020assessing}: the projected helical angle $\alpha$, that is the angle between the projection of the longitudinal fiber direction $\myVec{f}$ on the epicardial tangential plane defined by the circumferential and apico-basal directions $\refL$ and $\refN$; the intrusion angle $\gamma$, defined as the one between $\refL$ and the plane $\refL \text{-} \refN$; the cross-fiber angle $\beta$ defined as the rotation of the system around the $\refL$ direction (see Fig~\ref{fig:fig1}).

In the following paragraphs, we describe how the myocardial reference system is defined, how to compute the rotation angles from the experimental fiber reference system, and, vice-versa, how to reconstruct the fiber reference system from the angles.

\subsubsection{Myocardial reference system}\label{subsec:myo_RF}
\begin{figure}[t]
    \centering
    \includegraphics[width=\linewidth]{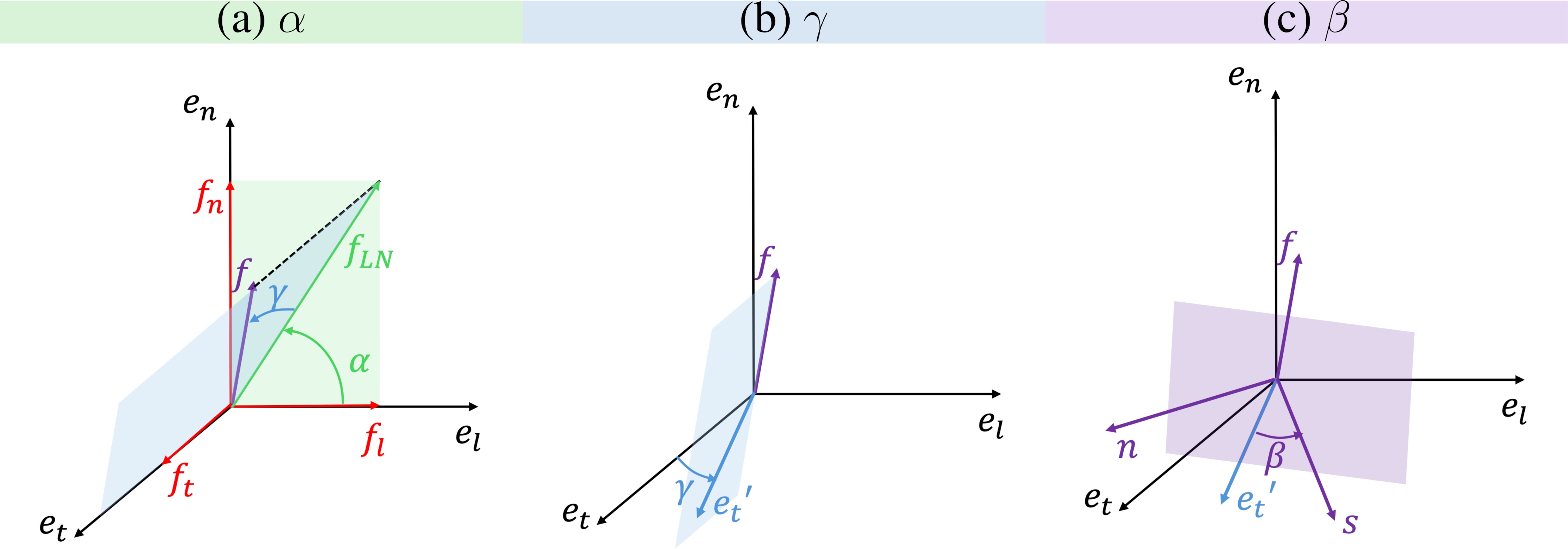}
    \caption{\textbf{Construction of the myocardial reference system.} (a) apico-basal distance $\psi$ (b) transmural distance $\phi$ , (c) myocardial reference system.}
    \label{fig:fig2}
\end{figure}

Several \acrshort{ldrbms} have been proposed in the literature, each employing different strategies to compute the myocardial reference system \cite{piersanti2021modeling}. In this work, we use the \acrshort{ldrbm} by Doste et al. \cite{doste2019rule, piersanti2021modeling}, where
the three directions $\refL$, $\refT$, and $\refN$ are computed from two scalar fields: the transmural distance $\phi$ defining the distance between the epicardium and the endocardial surfaces, see Fig~\ref{fig:fig2}(a); the apico-basal $\psi$ distance representing the distance from the apex towards the ventricular base, see Fig~\ref{fig:fig2}(b). Finally, the gradients of the transmural distance $\phi$ and the apico-basal distance $\psi$ are used to construct the transmural and apico-basal directions, respectively. 

The scalar fields $\phi$ and $\psi$ are computed by solving the following problems
\begin{align}
    \left\{
    \begin{aligned}
            &\Delta \phi = 0 &&\text{ in } \Omega\\
            &\phi = 0 &&\text{ on } \boardEpi, \\
            &\phi = -1 &&\text{ on } \boardEndoRV, \\
            &\phi = 2 &&\text{ on } \boardEndoLV,
    \end{aligned} 
    \right. &   &
    \left\{
    \begin{aligned}
            &\Delta \psi = 0 &&\text{ in } \Omega\\
            &\psi = 0 &&\text{ on } \boardApex, \\
            &\psi = 1 &&\text{ on } \boardBase,
    \end{aligned}
    \right.
\label{eq: phi_psi}
\end{align}
where $\Omega$ is the computational domain; $\boardEpi$, $\boardEndoRV$ and $\boardEndoLV$ represent the epicardium, 
the endocardium of the \acrfull{rv} and the endocardium of the \acrfull{lv}, respectively; while $\boardBase$ represents the ventricular base and $\boardApex$ the 
apex of the myocardium as shown in Fig S11. Given the variables $\phi$ and $\psi$, the transmural, apico-basal, and the circumferential directions $\myVec{e}_t$, $\refN$, and $\refL$ are computed according to
\begin{equation}
\begin{cases}
    \mathbf{e}_t = \dfrac{\nabla \phi}{\|\nabla \phi\|}, \\[4pt]
    \mathbf{e}_n = 
    \dfrac{
        \dfrac{\nabla \psi}{\|\nabla \psi\|}
        - \left( \dfrac{\nabla \psi}{\|\nabla \psi\|} \cdot \mathbf{e}_t \right)\mathbf{e}_t
    }{
        \left\|
        \dfrac{\nabla \psi}{\|\nabla \psi\|}
        - \left( \dfrac{\nabla \psi}{\|\nabla \psi\|} \cdot \mathbf{e}_t \right)\mathbf{e}_t
        \right\|
    }, \\[4pt]
    \mathbf{e}_l = \mathbf{e}_n \times \mathbf{e}_t.
\end{cases}
\label{eq:ref_sys}
\end{equation}

Figure~\ref{fig:fig2} illustrates the transmural distance $\phi$, the apico-basal distance $\psi$, and the myocardial reference system ($\refL$, $\myVec{e}_t$, $\refN$).

\subsubsection{From fibers to angles}
When both the fiber field and the myocardial reference system are available, the rotation angles $\alpha$, $\gamma$, and $\beta$ can be computed in sequence. To obtain $\alpha$, the longitudinal direction $\myVec{f}$ is projected  on the plane $\refN \text{-} \refL$, resulting in $\fLN$ as
\begin{equation}
\fLN = \myVec{f} - (\myVec{f} \cdot \myVec{e}_t) \myVec{e}_t,
\label{eq: f_ln}
\end{equation}
and than normalized in $\hat{\myVec{f}}_{LN}$
\begin{equation}
\hat{\myVec{f}}_{LN} = \frac{\fLN}{\| \fLN \|}.
\label{eq: norm_f_ln}
\end{equation}
After defining $\hat{\myVec{f}}_{LN}$, $\alpha$ is obtained as the angle between $\hat{\myVec{f}}_{LN}$ and $\refL$, with $\alpha \in (-\pi, \pi]$. This choice avoids treating angles $\alpha$ and $\alpha + 2\pi$ as distinct. Finally, the $\operatorname{atan2}$ function is used to compute $\alpha$
\begin{equation}
\alpha = \operatorname{atan2}(\refN \cdot \hat{\myVec{f}}_{LN}, \refL \cdot \hat{\myVec{f}}_{LN}),
\label{eq: alpha}
\end{equation}
where the sign of the dot product $\refL \cdot \hat{\myVec{f}}_{LN}$ specifies which half-plane contains $\hat{\myVec{f}}_{LN}$. Similarly, $\gamma$ is computed by inverting the sign for $|\alpha| > \pi / 2$:
\begin{equation}
    \gamma = \begin{cases}
        \operatorname{atan2}(\myVec{e}_t \cdot \myVec{f}, \myVec{f} \cdot \hat{\myVec{f}}_{LN}) &\text{ if } \lvert \alpha \rvert \leq \pi / 2, \\
        - \operatorname{atan2}(\myVec{e}_t \cdot \myVec{f}, \myVec{f} \cdot \hat{\myVec{f}}_{LN}) &\text{ if } \lvert \alpha \rvert > \pi / 2.
    \end{cases}
    \label{eq: gamma}
\end{equation}
To compute $\beta$, two auxiliary vectors are introduced: $\myVec{e}_t'$, which is the transmural $\myVec{e}_t$ direction rotated of $\gamma$ in the $\myVec{f} \text{-} \myVec{e}_t$ plane, and $\refN'$ which is the normal direction to the $\myVec{f} \text{-}\myVec{e}_t$ plane. Specifically, $\myVec{e}_t'$ and $\refN'$ are computed as
\begin{equation}
    \begin{cases}
        \displaystyle \myVec{e}_t' = \frac{\myVec{e}_t - (\myVec{e}_t \cdot \myVec{f})\myVec{f}}{\|  \myVec{e}_t - (\myVec{e}_t \cdot \myVec{f})\myVec{f} \|},\\[8pt]
        
        \displaystyle \refN' = \frac{\myVec{f} \times \myVec{e}_t}{\| \myVec{f} \times \myVec{e}_t \|}.
    \end{cases}
    \label{eq:auxiliary}
\end{equation}
Finally, the $\beta$ angle is evaluated as
\begin{equation}
     \beta = \begin{cases}
        \operatorname{atan2}(\myVec{e}_t' \cdot \myVec{s}, \myVec{s} \cdot \refN') &\text{ if } \alpha \geq 0, \\
        - \operatorname{atan2}(\myVec{e}_t' \cdot \myVec{s}, \myVec{s} \cdot \refN') &\text{ if } \alpha < 0. \end{cases}
        \label{eq:beta}
\end{equation}
Fig~\ref{fig:fig1} depicts the three rotations.

\subsubsection{From angles to fibers}\label{subsec:fib_to_angles}
When the angles are provided -- either measured experimentally from fiber orientations or computed by a \acrshort{ldrbm} -- the myofiber field is reconstructed via the inverse process, which converts the angles back to the fiber reference system. To compute the fiber reference system $\myVec{f}$, $\myVec{s}$, $\myVec{n}$, from the values of $\alpha$ and $\gamma$ two rotations are imposed: the first, around direction $\myVec{e}_t$, is used to compute $\fLN$; the second, on the $\fLN \text{-} \myVec{s}$ plane, maps $\fLN$ into $\myVec{f}$ and $\myVec{e}_t$ into $\myVec{e}_t'$. In addition, $\refNPrime$ is defined as the cross product of $\refTPrime$ and $\myVec{f}$ as follows
\begin{equation}
\begin{cases}
    \myVec{f}
    = \cos(\alpha) \cos(\gamma) \, \refL 
    + \sin(\alpha) \cos(\gamma) \, \refN 
    + \sin(\gamma) \, \myVec{e}_t, \\[4pt]
    
    \refTPrime
    = \dfrac{
        \myVec{e}_t - \sin(\gamma) \, \myVec{f}
    }{
        \left\| 
            \mathbf{e}_t - \sin(\gamma) \, \myVec{f} 
        \right\|
    }, \\[8pt]
    
    \refNPrime
    = \refTPrime\times \myVec{f}.
\end{cases}
\label{eq:recSys}
\end{equation}
If $\beta = 0$, then $\myVec{s} = \myVec{e}_t'$ and $\myVec{n} = \refN'$, otherwise if $\beta \neq 0$ the rotation around the longitudinal direction $\myVec{f}$ have to be considered by multiplying the fiber reference system with the rotation matrix $\myTens{R}_{\beta}$ as follows
\begin{equation}
    \begin{bmatrix} 
        \myVec{f} \\
        \myVec{e}_t'\\
        \refN'
    \end{bmatrix} = 
    \myTens{R}_{\beta}
    \begin{bmatrix}
        \myVec{f} \\
        \myVec{s}\\
        \myVec{n}
    \end{bmatrix},
\end{equation}
where $\myTens{R}_{\beta}$ is defined as
\begin{equation}
\myTens{R}_{\beta} = \begin{bmatrix}
             &1 &0 &0 \\
             &0 &\cos(\beta) &-\sin(\beta)\\
             &0 &\sin(\beta)  &\cos(\beta)
        \end{bmatrix}  .
        \label{eq:rot_mat}
\end{equation}

 \subsection{Fibers measurement}\label{subsec:fib_measurement}
The data used in this work have been measured by the optical imaging technique proposed in \cite{giardini2025correlative}.  This approach provides high-resolution information on the three-dimensional organization of the myocardium, including the longitudinal and sheet directions of cardiomyocytes, across the entire mouse heart. Briefly, whole mouse hearts were fixed and optically cleared using the SHIELD protocol optimized for cardiac tissue \cite{olianti2022optical}. Cleared hearts were imaged using a modified mesoscopic selective plane illumination microscope (mesoSPIM), exploiting the intrinsic autofluorescence of cardiac muscle to reconstruct the myocardial architecture at near-isotropic spatial resolution ($\SI{3.25}{\micro\meter} \times \SI{3.25}{\micro\meter} \times \SI{3}{\micro\meter}$) \cite{giardini2025correlative}. 
The resulting volumetric datasets allowed a detailed three-dimensional reconstruction of the myocardium at single-cell scale.  Local cardiomyocyte orientation was then quantified with an isotropic spatial resolution of $\SI{96}{\micro\meter}$ by applying the Structure Tensor Analysis (STA) to the autofluorescence signal: the longitudinal myofiber field and the sheet direction were obtained from the eigenvectors associated with the smallest and the second smallest eigenvalues, respectively. The resulting STA-based mapping represents myocardial organization through two unit-length vector fields, describing the longitudinal myofiber field and the sheet orientation field.

\subsection{Decoupling macroscopic fiber orientations from microscopic disarray}\label{subsec:decoupling}
In this section, we present the methodology employed to derive the main fiber organization from an experimentally measured myofiber field, by separating macroscopic architecture from microscopic fiber disarray. The underlying idea is to perform a decoupling in the spatial frequency domain, treating high-frequency components as indicative of fiber disarray, while associating the remaining low-frequency content with the macroscopic fiber architecture. To this end, a Helmholtz filter is applied to the fiber field. This filter is particularly well suited for data defined on unstructured meshes and irregular domains, as encountered in cardiac geometries. Acting as a low-pass spatial filter, the Helmholtz filter attenuates high-frequency variations in fiber orientation, thereby removing fiber disarray. The resulting smoothed field is interpreted as the macroscopic fiber architecture, while the fiber disarray is subsequently defined as the residual signal.

This procedure can be applied independently to the three angles $\alpha$, $\beta$, and $\gamma$. 
For the sake of simplicity, in what follows we focus on the $\alpha$ angle. The measured angle field $\alpha(\myVec{x})$ is decomposed into a macroscopic component $\tilde{\alpha}(\myVec{x},\ell)$ and a microscopic disarray term $\varepsilon_{\alpha}(\myVec{x},\ell)$:
\begin{equation}
\alpha(\myVec{x}) = \tilde{\alpha}(\myVec{x}, \ell) + \varepsilon_{\alpha}(\myVec{x}, \ell),
\label{eq:decouple}
\end{equation}
where, $\ell$ denotes the regularization radius representing the characteristic length scale of the smoothing operator used to extract the macroscopic field.

A critical issue is that direct smoothing in the angular domain is not straightforward for two primary reasons. First, angular variables are circular, meaning that $\alpha$ and $\alpha + 2\pi$ represent the same value. Second, due to the directional invariance of the fibers, angles differing by $\pi$ represent the same physical fiber direction, and therefore $\alpha$ and $\alpha + \pi$ should be treated as equivalent.
For these reasons, smoothing is not performed directly on $\alpha$, but rather in a transformed coordinate system that correctly reflects this topology. Specifically, the fiber angles are mapped onto the coordinate system $(\eta_{\sin}, \eta_{\cos})$, defined as
\begin{equation}
\begin{cases}
    \eta_{\sin} = \sin(2\alpha), \\
    \eta_{\cos} = \cos(2\alpha).
\end{cases}
\label{eq: cartesian}
\end{equation}
Since $\sin(2\alpha) = \sin\!\left(2(\alpha+\pi)\right)$ and $\cos(2\alpha) = \cos\!\left(2(\alpha+\pi)\right)$, opposite fiber directions are mapped to the same point in this coordinate system.

By applying a Helmholtz filter with regularization radius $\ell$, the regularized coordinates $\tilde{\eta}_{\sin}$ and $\tilde{\eta}_{\cos}$ are obtained by solving the problems:
\begin{equation}
\begin{cases}
    -\ell^2 \Delta \tilde{\eta}_{\sin} + \tilde{\eta}_{\sin} = \eta_{\sin} & \text{ in } \Omega,\\
    \nabla \tilde{\eta}_{\sin} \cdot \mathrm{n}_{\Gamma} = 0 & \text{ on } \Gamma,
\end{cases}
\text{  and   }
\begin{cases}
    -\ell^2 \Delta \tilde{\eta}_{\cos} + \tilde{\eta}_{\cos} = \eta_{\cos}  & \text{ in } \Omega, \\
    \nabla \tilde{\eta}_{\cos} \cdot \mathrm{n}_{\Gamma} = 0 & \text{ on } \Gamma,
\end{cases}
\label{eq: regularization}
\end{equation}
where $\Gamma$ is the whole boundary of the biventricular domain $\Omega$. In an unbounded domain, each problem of Eq.~\eqref{eq: regularization} is equivalent to applying the convolution operator
\begin{equation}
     \tilde{\eta}(\myVec{x}) = (\eta * h)(\myVec{x}) 
    = \int_{\Omega} h\!\left(\lvert \myVec{x} - \myVec{x}' \rvert\right) 
    \eta(\myVec{x}') \, \mathrm{d}\myVec{x}',
\label{eq:conv}
\end{equation}
where $h(r)$ is the Green’s function of the operator $I - \ell^2 \Delta$, given by
\begin{equation}
    h(r) = \frac{1}{4 \pi \ell^2} 
    \frac{e^{-r / \ell}}{r}.
\label{eq:kernel}
\end{equation}
As shown in Eq.~\eqref{eq:kernel}, the parameter $\ell$ defines the characteristic smoothing length and controls the attenuation of spatial frequencies significantly higher than $1/\ell$, effectively acting as a low-pass spatial filter.
In practice, Eq.~\eqref{eq: regularization} is solved numerically using the finite element method, which is well suited for unstructured meshes and irregular computational domains.

Finally, the smoothed angle field $\tilde{\alpha}$ (in radians) is recovered by mapping the regularized coordinates back to the angular domain using the $\operatorname{atan2}$ operator:
\begin{equation}
\tilde{\alpha} = \frac{1}{2}\operatorname{atan2}(\tilde{\eta}_{\sin}, \tilde{\eta}_{\cos}).
\label{eq: smoothed_alpha}
\end{equation}
After obtaining the regularized angle field, the fiber disarray $\varepsilon_{\alpha}$ is determined directly from Eq.~\eqref{eq:decouple}.

As mentioned, the same procedure can be applied to all three angles $\alpha$, $\beta$, and $\gamma$. 
However, in the present work, due to the transversely isotropic behavior of the material, rotations in the $\myVec{s}$--$\myVec{n}$ plane, represented by the angle $\beta$, are neglected. 
Therefore, after computing $\tilde{\alpha}$ and $\tilde{\gamma}$, the regularized fiber reference system $\tilde{\myVec{f}}$, $\tilde{\myVec{s}}$, and $\tilde{\myVec{n}}$ is reconstructed using Eq.~\eqref{eq:recSys}.

\begin{figure}[t]
\centering
\includegraphics[width=\linewidth]{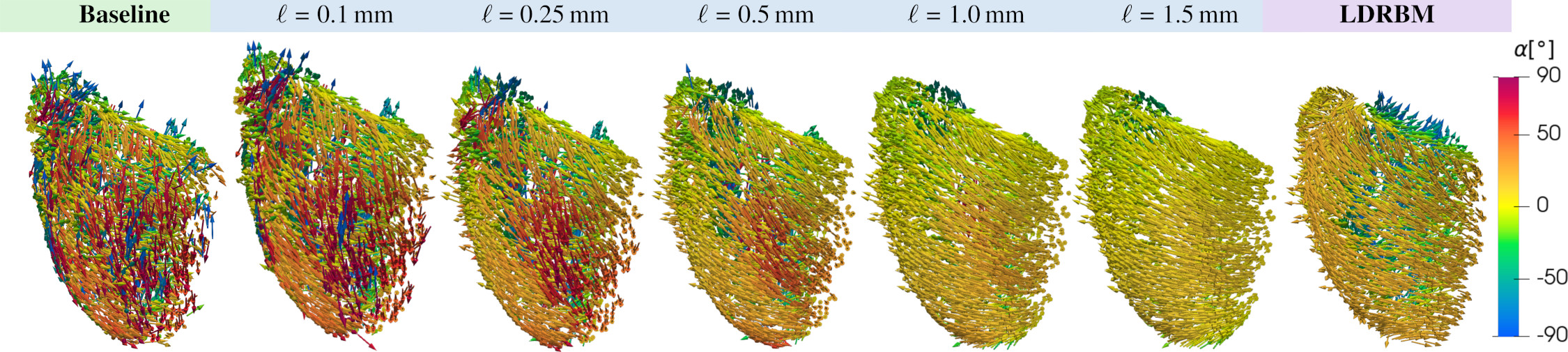}
\caption{\textbf{Comparison different fiber fields derived by varying the regularization radius $\ell$ with respect to the the unfiltered fiber field and the \acrshort{ldrbm} one.} On the left the unfiltered fiber field and on the right the \acrshort{ldrbm} one \cite{doste2019rule}. The colormap represents the angular variation of $\alpha$ within the myocardium.}
 \label{fig:fig3}
\end{figure}

Figure~\ref{fig:fig3} compares the baseline fiber field obtained from experimental measurements with the corresponding regularized fields for increasing values of the regularization radius $\ell$, as well as with the fiber field generated by the \acrshort{ldrbm}. In the present work, the maximum diameter of the base is of approximately $\SI{8}{\milli \meter}$.

\subsection{Prescribing the fiber orientations via \acrshort{ldrbm}}
The Doste et al. \acrshort{rbm} \cite{piersanti2021modeling, doste2019rule} surrogates the fiber architecture by interpolating specific angular values of $\alpha$ and $\gamma$ at the epicardium and endocardium and linearly interpolating them along the transmural direction. A convex combination is performed depending on $\phi$ as follows
\begin{equation}
\alphaRBM =
\begin{cases}
    -\,\alphaEndoRV \, \phi 
    + \alphaEpiRV \, (\phi - 1), 
    & \text{in } \Omega_{RV}, \\[6pt]
    \alphaEndoLV \, \phi 
    + \dfrac{\alphaEpiLV}{2} \, (1 - \phi), 
    & \text{in } \Omega_{LV},
\end{cases}
\label{eq:conv_combination}
\end{equation}
where $\alpha_{endo, RV}$, $\alpha_{endo, LV}$, $\alpha_{epi, LV}$ and $\alpha_{epi, RV}$ are the prescribed angles at right and left endocarium and epicardium, respectively, while $\Omega_{LV}$  and $\Omega_{RV}$ denote the left and right ventricular domains. This procedure is applied to both the $\alpha$ and $\gamma$ angles. The resulting angles  $\alphaRBM$ and $\gammaRBM$ are then used in Eq.~\eqref{eq:recSys} to construct the myocardial reference system.

To identify $\alphaEndoRV$, $\alphaEndoLV$, $\alphaEpiRV$, $\alphaEpiLV$, $\gammaEndoRV$, $\gammaEndoLV$, $\gammaEpiRV$ and $\gammaEpiLV$ that best match an experimentally measured fiber field, we compute the statistical distribution of $\alpha$ and $\gamma$ selecting the modal value in different myocardial sub-regions, as done in \cite{piersanti2025defining} (see Fig \ref{fig:histograms}).

\begin{figure}[t]
    \centering
    \includegraphics[width=\linewidth]{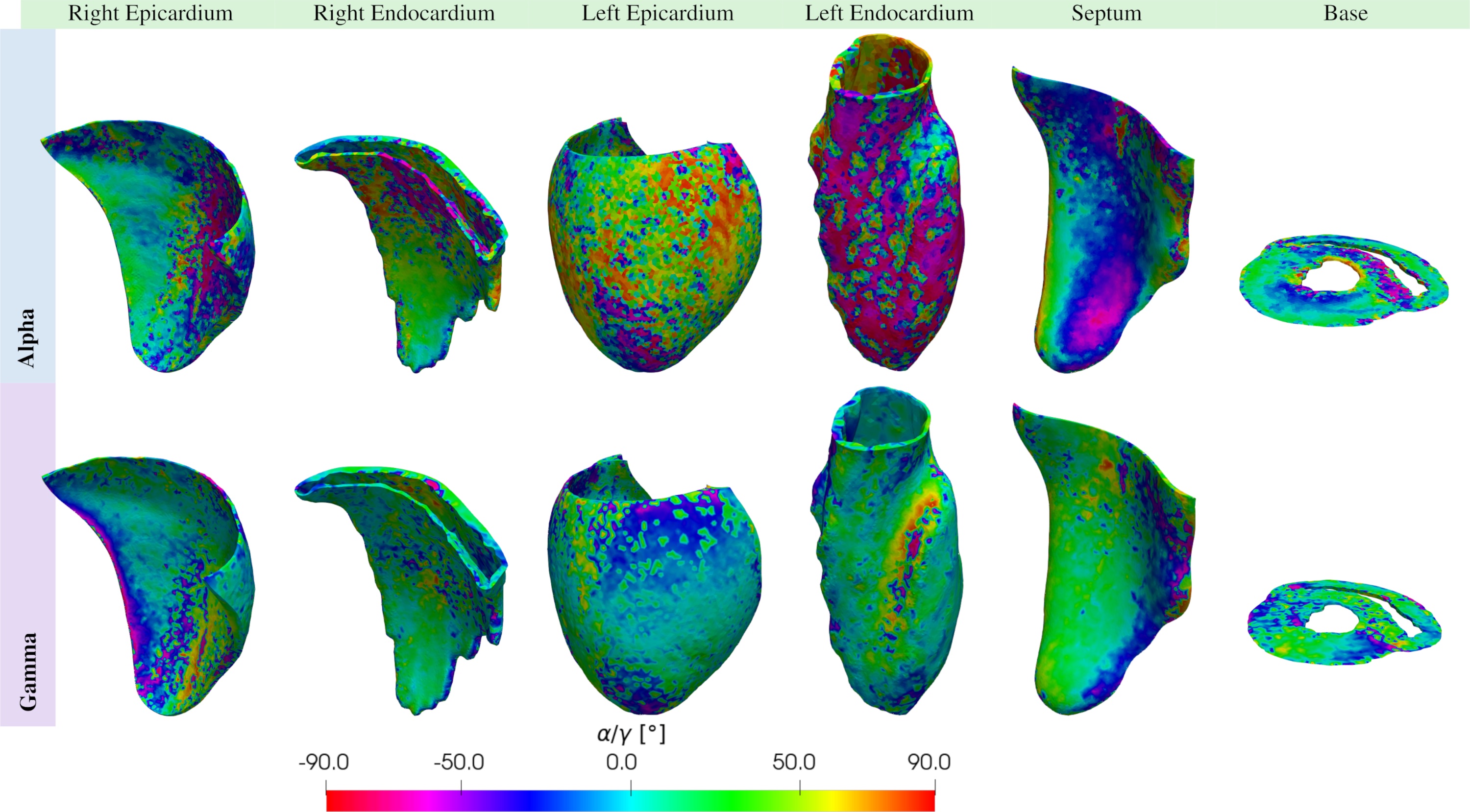}
    \caption{\textbf{Angles mapping.} Variation of $\alpha$ and $\gamma$ angles within the ventricular sub-regions.}
    \label{fig:alphaMap}
\end{figure}

\begin{figure}[t]
    \centering
    \includegraphics[width=\linewidth]{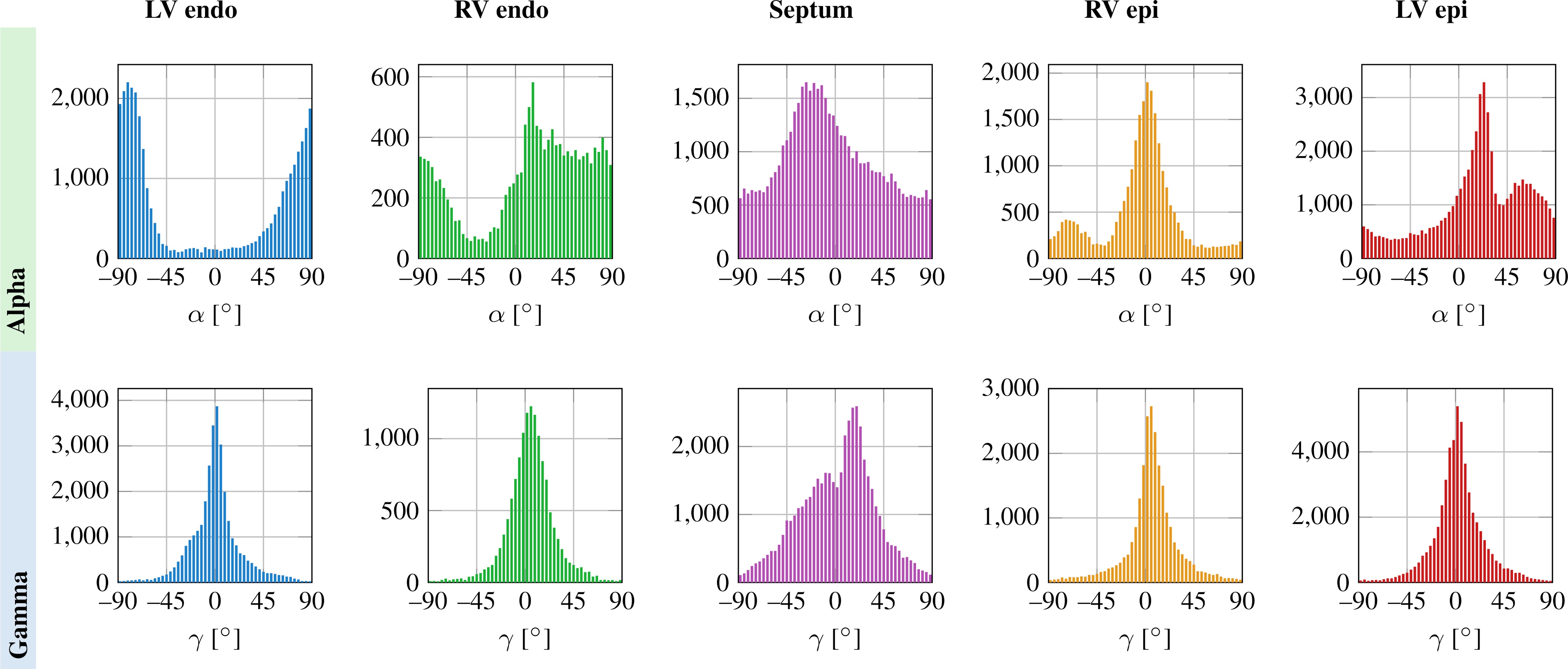}
    \caption{\textbf{Angles distribution.} Distribution of the $\alpha$ and $\gamma$ angles in different zones of the domain.}
    \label{fig:histograms}
\end{figure}

Figure~\ref{fig:alphaMap} reports the values of the angles $\alpha$ and $\gamma$ across different regions of the ventricular domain. Compared with $\alpha$, $\gamma$ is relatively uniform and remains near zero. This suggests that global ventricular contraction is primarily driven by rotations in the $\refL$–$\refN$ plane, which exhibit a pronounced transmural variation from endocardium to epicardium.
Moreover, $\alpha$ displays significant variability even within the same region, suggesting that assigning a single fiber angle value within the entire epicardium and endocardium may not be sufficient to accurately represent local myocardial fiber orientation. Finally, it is worth noting that the \acrshort{ldrbm} by Doste et al. \cite{doste2019rule} does not provide a reconstruction of the $\beta$ angle, assuming a transversely isotropic material behavior.

\section{Electromechanical model}\label{sec:Baseline}
The electromechanical model consists of multiple core modules, each representing different physical processes: electrophysiology, activation, mechanics, and blood circulation. The electrophysiology component determines the activation time, defined as the instant at which the electrical stimulus reaches each point in the cardiac domain, as well as the calcium ion concentration. These quantities are then used as input to the activation core model, which computes the active tension along the fibers by accounting for both calcium concentration and local deformations. The mechanics combines the contributions of active tension and passive tissue behavior to compute displacements, while the blood circulation component, which represents the hemodynamics of the entire cardiovascular system, provides the pressure boundary conditions at the endocardium \cite{fedele2023comprehensive,regazzoni2022cardiac,piersanti20223d}.

Under the physiological conditions considered, the electrophysiology component is described using an eikonal-diffusion formulation \cite{franzone2014mathematical}. This enables the efficient computation of activation times by solving an eikonal equation. The local calcium concentration is computed by evaluating an experimentally derived calcium transient at the activation time associated with each point, following the approach proposed in \cite{stella2022fast}. The electrical signal propagation across the myocardial domain is governed by the conductivity tensor $\myTens{D}$, defined as
\begin{equation}
    \myTens{D} = \sigmaF\frac{\myTens{F} \myVec{f}_0 \otimes \myTens{F} \myVec{f}_0}{\left| \left| \myTens{F} \myVec{f}_0 \right| \right|^2} +
    \sigmaS\frac{\myTens{F} \myVec{s}_0 \otimes \myTens{F} \myVec{s}_0}{\left| \left| \myTens{F} \myVec{s}_0 \right| \right|^2} + 
    \sigmaN\frac{\myTens{F} \myVec{n}_0 \otimes \myTens{F} \myVec{n}_0}{\left| \left| \myTens{F} \myVec{n}_0 \right| \right|^2},
    \label{eq:conductvity}
\end{equation}
where $\sigmaF$, $\sigmaS$, and $\sigmaN$ are the conductivities in the fiber, sheet, and normal directions, respectively. In modeling practice, $\sigmaF$ is typically assumed to be twice as large as $\sigma_S$, with $\sigma_S$ set equal to $\sigma_N$. Hence, the fiber orientation $\myVec{f}_0$ establishes a preferential direction for signal propagation. 

The cardiac mechanics, describing the dynamics of the tissue displacement $\myVec{d}$, is modeled using the momentum conservation equation under the hyperelasticity assumption, combined with an active stress approach \cite{guccione1991finite,fedele2023comprehensive}.
The deformation gradient tensor, computed from the displacement field $\myVec{d}(\myVec{x})$, is defined as $\myTens{F} = \myTens{I} + \nabla \myVec{d}(\myVec{x})$.

The constitutive behavior of the myocardium is governed by the Piola-Kirchhoff stress tensor, defined as
\begin{equation}
    \myTens{P} = \frac{\partial \mathcal{W}(\myTens{F})}{\partial \myTens{F}} 
    + \mF \ActTension\frac{\myTens{F}\myVec{f}_0 \otimes \myVec{f}_0}{\lvert \lvert \myTens{F}\myVec{f}_0\rvert \rvert}
    + \mS \ActTension\frac{\myTens{F}\myVec{s}_0 \otimes \myVec{s}_0}{\lvert \lvert \myTens{F}\myVec{s}_0\rvert \rvert} 
    + \mN \ActTension\frac{\myTens{F}\myVec{n}_0 \otimes \myVec{n}_0}{\lvert \lvert \myTens{F}\myVec{n}_0\rvert \rvert},
    \label{eq:Piola}
\end{equation}
where $\mathcal{W}$ is the Usyk strain potential energy established with transversely isotropic parameters \cite{usyk2000effect}; $\ActTension$ is the active tension, and $\mF$, $\mS$, and $\mN$ are the Stress Factors ($\SF$) in the fiber, sheet, and normal directions, respectively \cite{piersanti20223d}. 
The $\SF$ configuration $\mF = 1$, $\mS = \mN = 0$ corresponds to an active contraction occurring solely along the fiber longitudinal direction.
In computational modeling practice, the other components can be activated to surrogate the effect of fiber disarray \cite{guan2020effect,piersanti20223d}.

The activation provided by the RDQ20-MF model proposed by \cite{regazzoni2020biophysically}. The latter computes the active tension field by solving a system of 20 ordinary differential equations. 
The RDQ20-MF model offers accurate results under physiological conditions at low computational cost, while accurately capturing cooperative, length-dependent activation and force-velocity relationships \cite{regazzoni2020biophysically,fedele2023comprehensive}. The output of the RDQ20-MF model is the active tension $\ActTension$, defined as
\begin{equation}
    \ActTension = \CBstiffness \, \mathcal{G}(\Calcium, \SL),
    \label{eq:Ta}
\end{equation}
where $\CBstiffness$ is the crossbridge stiffness, which linearly scales the active tension with respect to the crossbridge contraction in the sarcomere, $\mathcal{G}$; 
$\Calcium_i$ is the calcium ion concentration provided by the electrophysiology model; 
and $\SL$ is the sarcomere length, resulting from fiber contraction computed in the mechanics.
The crossbridge stiffness $\CBstiffness$ is set to $\SI{23}{\mega\pascal}$ baseline simulation. This parameter can be adjusted to reduce or increase the global contractility of the myocardium. 
Furthermore, the $\CBstiffness$ in the \acrshort{rv} is half of that in the \acrshort{lv}, effectively imposing a contractility ratio of $1/2$ between \acrshort{lv} and \acrshort{rv}.

Mechanical boundary conditions are imposed in four cardiac zones: epicardium, left and right endocardium, and the ventricular base. Robin-type boundary conditions are applied at the epicardium and basal plane, following the approach in \cite{pfaller2019importance}.  
At the endocardium, the mechanical model is coupled with the blood circulation model through energy-consistent boundary conditions, 
as proposed in \cite{regazzoni2022cardiac}.

The blood circulation of the entire cardiovascular system is modeled using a 0D closed-loop approach, as proposed in \cite{regazzoni2022cardiac}. The systemic and pulmonary circulations are represented by 0D Resistance–Inductance–Capacitance (RLC) networks, in which blood flow rate is analogous to electrical current, and blood volume corresponds to electrical potential. Heart chambers are modeled as time-varying elastance elements, while non-ideal diodes represent the heart valves \cite{regazzoni2022cardiac,piersanti20223d,fedele2023comprehensive}.

\subsection{Numerical Implementation}
The numerical framework has been implemented within \texttt{life}${}^\text{\texttt{x}}$ (\url{https://lifex.gitlab.io})\cite{africa2022lifex,africa2023lifexep,africa2023lifexfiber,bucelli2025lifex}, an in-house high-performance C++ \acrfull{fe} library focused on
cardiac applications based on \texttt{deal.II} FE core (\url{https://www.dealii.org})\cite{arndt2017deal,arndt2021deal}. The simulations were run on one CPU node powered by two 24-core CPUs with 512 GB RAM.

The multiphysics problem is discretized on a tetrahedral mesh with an element size of $\SI{150}{\micro\meter}$ and composed of 328000 nodes.  
Space discretization is based on linear \acrshort{fe}, and the time discretization is performed using a Backward Difference Formula (BDF) of second order for the acceleration term. We solve the non-linear equations at each timestep using the Newton method. Finally, for the blood circulation, a forward Euler scheme is employed. The circulation model is then coupled with the mechanics through a Lagrange multiplier formulation, in which the pressures of the \acrshort{lv} and \acrshort{rv} act as Lagrange multipliers in the coupled circulation-mechanics problem \cite{regazzoni2022cardiac,piersanti20223d,fedele2023comprehensive}.

The timestep for electromechanics and blood circulation is set to $\SI{1}{\milli\second}$. Each heartbeat for the physiological mouse heart lasts $\SI{200}{\milli\second}$, and $5$ heartbeats are simulated. The ectrophysiological parameters were calibrated using a heart-specific activation map obtained from experimental measurements, and the activation model was tuned based on shortening twitch test results. To calibrate the remaining parameters, we relied on physiological data of blood circulation \acrfull{qoi} derived from the literature \cite{tabima2010measuring,townsend2016measuring}. 
Furthermore, to reduce the computational time and resources required during the tuning phase, we employed the 0D emulator proposed by \cite{regazzoni2021accelerating}. 
This emulator allows the computation of convenient initial conditions for the 3D simulation, enabling faster convergence toward the limit cycle \cite{piersanti20223d}. It can also be used to calibrate the parameters of the RLC circuit representing blood circulation, as well as selected parameters of the 3D model. Further details are provided in Appendix~\ref{app:Baseline_sim}.
\section{Results}\label{sec:Results}

In this section, we report the numerical results obtained with the cardiac electromechanical model, focusing on the role of the fiber field and the associated fiber disarray on the electrophysiology (Sec.~\ref{subsec:electrophysiology}), passive mechanics (Sec.~\ref{subsec:passive_mec}), and finally on the overall electromechanical function (\ref{subsec:electromechanics}). In this analysis, the baseline simulation, built on the experimental fiber field, is taken as reference result and simulations built on regularized or \acrshort{ldrbm} fiber fields are compared with it.

We remark that, in the murine heart, global mechanical function is predominantly governed by the \acrshort{lv}, whereas the \acrshort{rv} operates under substantially lower pressure and mechanical load as a consequence of the low-resistance pulmonary circulation \cite{milani2014small}. 
Moreover, owing to its extremely thin free wall and strong ventricular interdependence, \acrshort{rv} dynamics are largely driven by \acrshort{lv} contraction.
From an experimental standpoint, these features also lead to a reduced signal-to-noise ratio in \acrshort{rv} fiber measurements. For these reasons, and given the secondary contribution of the \acrshort{rv} to overall cardiac mechanics in murine hearts, the present analysis focuses primarily on the \acrshort{lv}.

\begin{figure}[t]
    \centering
    \includegraphics[width=\linewidth]{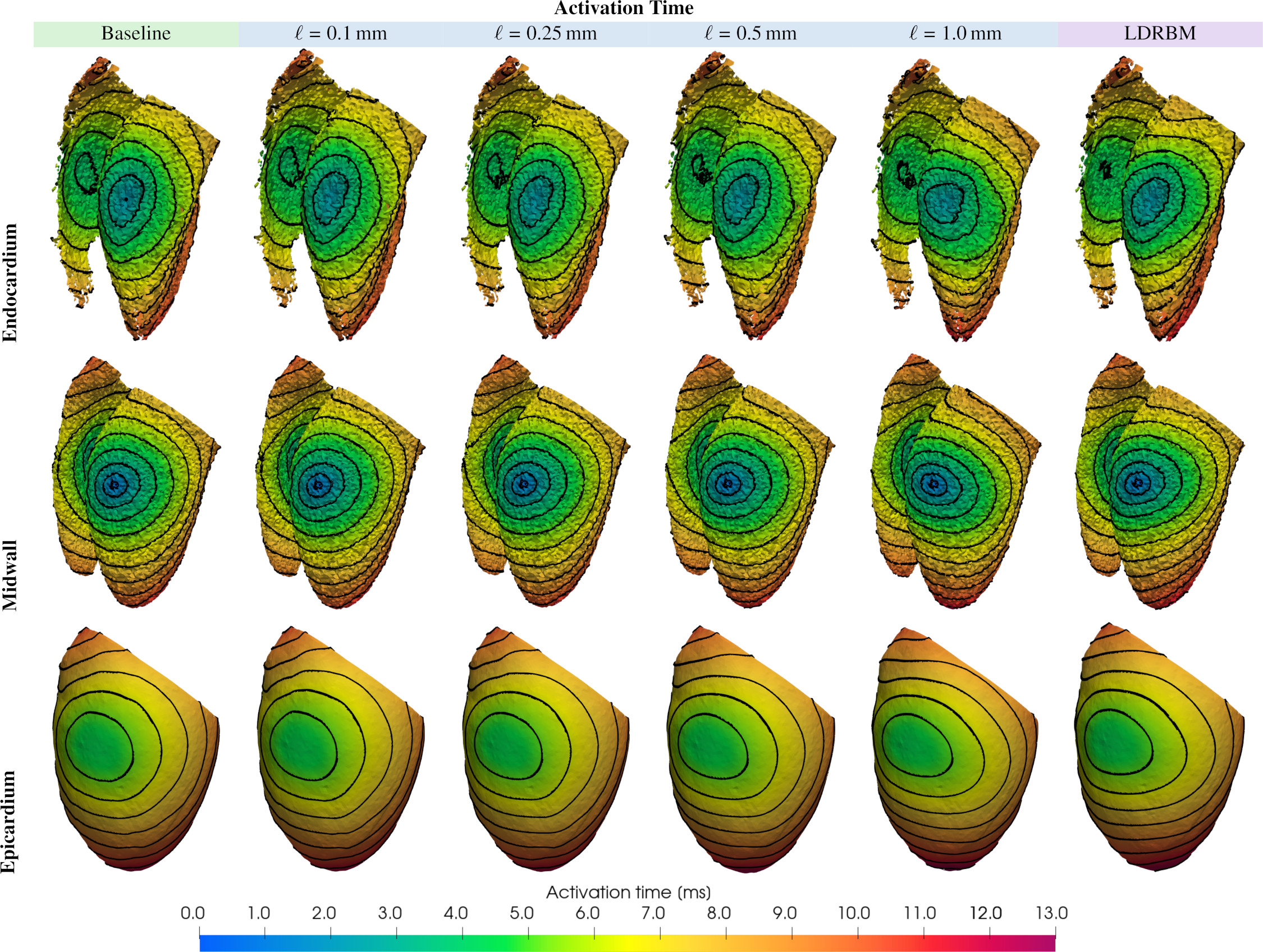}
    \caption{\textbf{Influence of the fiber field on electrophysiology propagation in the myocardium.} The distance between the isochrones is \SI{1.0}{\milli \second} Endocardial and epicardial activation maps are shown for different fiber architectures: experimental (baseline), varying regularization radii $\ell$, and the \acrshort{ldrbm} fibers.}
    \label{fig:ATMap}
\end{figure}

\subsection{Electrophysiology}\label{subsec:electrophysiology}
We evaluated the influence of the fiber field on the propagation of the electrophysiology signal through the myocardium. Tissue conductivity was modeled as transversely isotropic, with fiber-longitudinal conductivity twice that of the sheet and normal directions.

Fig~\ref{fig:ATMap} shows the endocardial and epicardial activation maps for increasing fiber regularization radius and for the \acrshort{ldrbm}. The results for the baseline case are aligned with the literature \cite{crocini2016optogenetics}. Wavefronts originating from the stimulation sites propagate preferentially along the higher-conductivity fiber direction. As expected, in the experimental fiber field, the presence of disarray reduced the directional effect of faster conduction along the fibers, resulting in a more isotropic propagation pattern compared to \acrshort{ldrbm} fibers. As expected, with the experimental fiber field, the presence of disarray mitigated the effect of faster conduction in the fiber direction. Moreover, as the fiber regularization radius increases, the preferred direction of propagation in the longitudinal direction becomes more pronounced. Overall, the \acrshort{ldrbm} produces an activation pattern similar to that of the regularized fiber field, suggesting that \acrshort{ldrbm} captures the mean fiber orientation. However, signal propagation with \acrshort{ldrbm} is more anisotropic compared to the experimental fiber field. Despite increased regularization of the fiber field enhances anisotropic conductivity, the latest activation times (around $\SI{13}{\milli \second}$) show no clear dependence on fiber smoothing or on the rule-based reconstruction.

\subsection{Passive Mechanics}\label{subsec:passive_mec}
We investigated how the fiber architecture affects passive inflation of the biventricular geometry. The following conditions were considered: no active stress, prescribed cavity pressures without explicitly modeling the blood circulation. 
The resulting numerical problem was solved by Newton-Raphson iterations using pressure as the control variable. Pressure was progressively applied to both the \acrshort{lv} and \acrshort{rv}, increasing from $0$ to $\SI{37.5}{\mmHg}$ in $50$ equal increments. This approach allowed us to examine the cardiac passive response across a range that includes physiological diastolic pressures.

\begin{figure}[t]
    \centering
    \includegraphics[width=\linewidth]{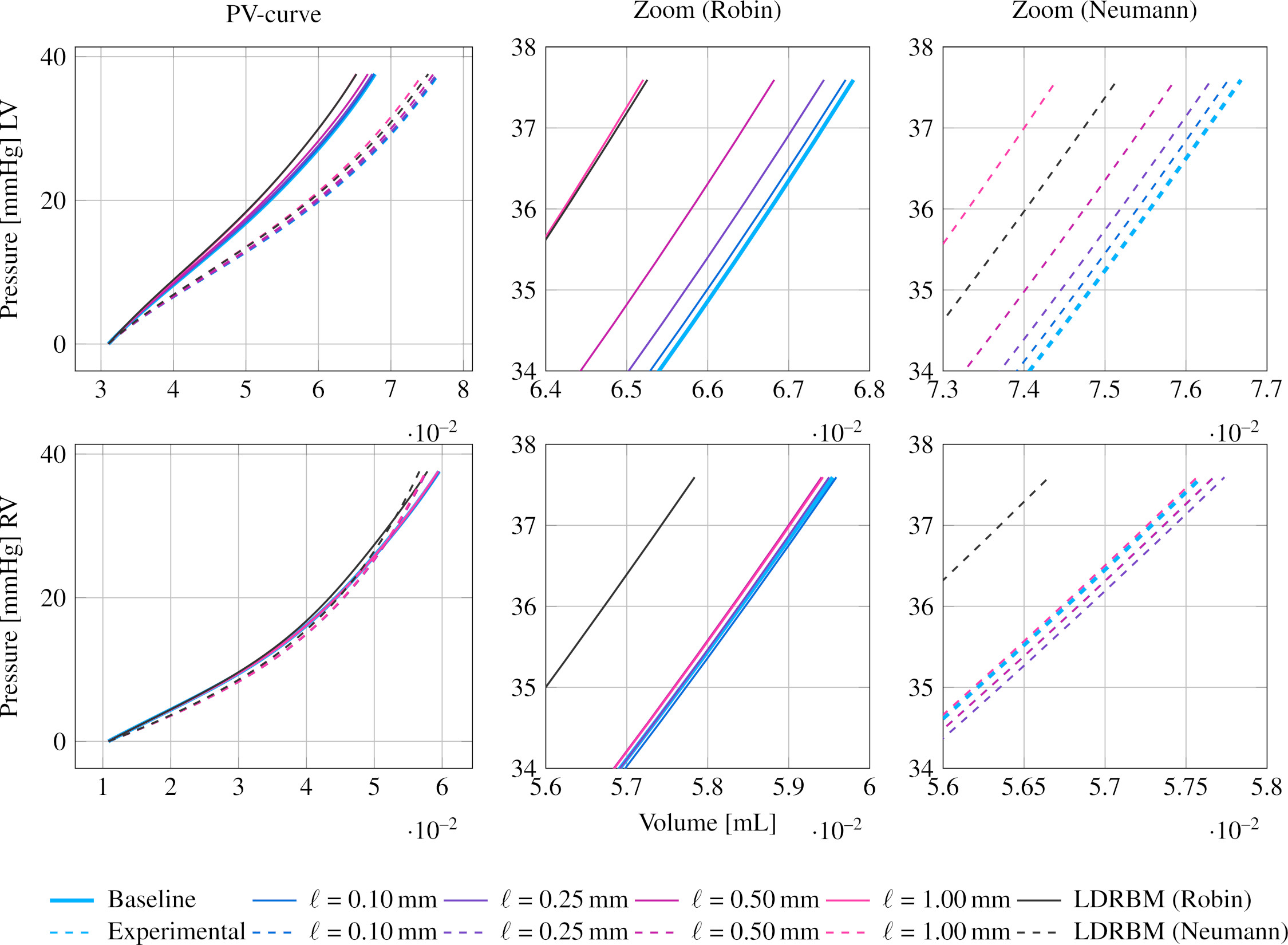}
    \caption{\textbf{Passive inflation of the left (PV-curve LV) and right (PV-curve RV) ventricles under increasing pressure, in absence of active contraction.} Comparison among different fiber architectures: experimental (baseline), progressively smoothed fibers with increasing regularization radius $\ell$, and \acrshort{ldrbm}. Two cases are shown. The first one (solid lines) correspond to the baseline \acrshort{bc} with Robin \acrshort{bc} at the base and at the epicardium. The second case (dashed lines) corresponds to homogeneous Dirichlet \acrshort{bc} at the base and Neumann \acrshort{bc} at the epicardium. The second column shows the zoomed views of the Robin \acrshort{bc} case, the third column shows the zoomed views of the Neumann \acrshort{bc} case.}
    \label{fig:PassMec}
\end{figure}

Fig \ref{fig:PassMec} compares the inflation curves of \acrshort{lv} and \acrshort{rv} obtained with different fiber configurations: 
experimental, progressively smoothed fields (with increasing regularization radius, as described in \ref{subsec:decoupling}, see also Fig~\ref{fig:fig3}), and Doste el al. \acrshort{ldrbm} \cite{doste2019rule}. Moreover, we analyzed two different cases. In the first case, we used Robin \acrshort{bc} at the base and at the epicardium, as in the baseline simulation (see Tab.\ref{tab:TabBaseline}); in the second case we replaced the Robin \acrshort{bc} at the epicardium with stress-free (Neumann) \acrshort{bc}, and applied homogeneous Dirichlet \acrshort{bc} at the base to constraint the motion. 
The rationale behind the letter case was to reduce the effect of \acrshort{bc} on the elastances of the chambers, thus isolating the effect of the fiber field. 
The results show that, in both the cases considered, reduced fiber dispersion (i.e., greater alignment coherence) slightly increases myocardial stiffness. More precisely, increasing $\ell$ from $\SI{0}{\milli \meter}$ to $\ell = \SI{1.0}{\milli \meter}$ reduces \acrshort{lv} volume variation at $P = \SI{37.5}{\mega \pascal}$ by 4\% with Robin \acrshort{bc} and by 3\% with Neumann \acrshort{bc}. 
Moreover, replacing the experimental fiber field with the \acrshort{ldrbm} one yields a similar variations.
Since these trends are consistent across boundary-condition setups, we conclude that the observed differences can be attributed primarily to fiber architecture and disarray, rather than to the influence of the surrounding tissues.

\subsection{Electromechanics}\label{subsec:electromechanics}
We analyzed the effects of the fiber architecture on the overall electromechanical behavior. We stress that when electromechanical function is considered, the fiber architecture influences all the major components of the heart model. Specifically, fiber orientation affects passive mechanics and electrophysiology, as previously discussed (in Sec. 4.1 and 4.2), and also plays a pivotal role in the active stress (see Sec. 3 and Eq. 20).

To assess the impact of fibers on the electromechanical function, we employ two types of readout: pressure-volume data -- providing global indication of the pumping function of the heart -- and strain maps -- yielding local information on the tissue mechanical response.

\subsubsection{Pressure-Volume analysis}\label{subsec:pv_analysis}

\begin{figure}[t]
    \centering
    \includegraphics[width=\linewidth]{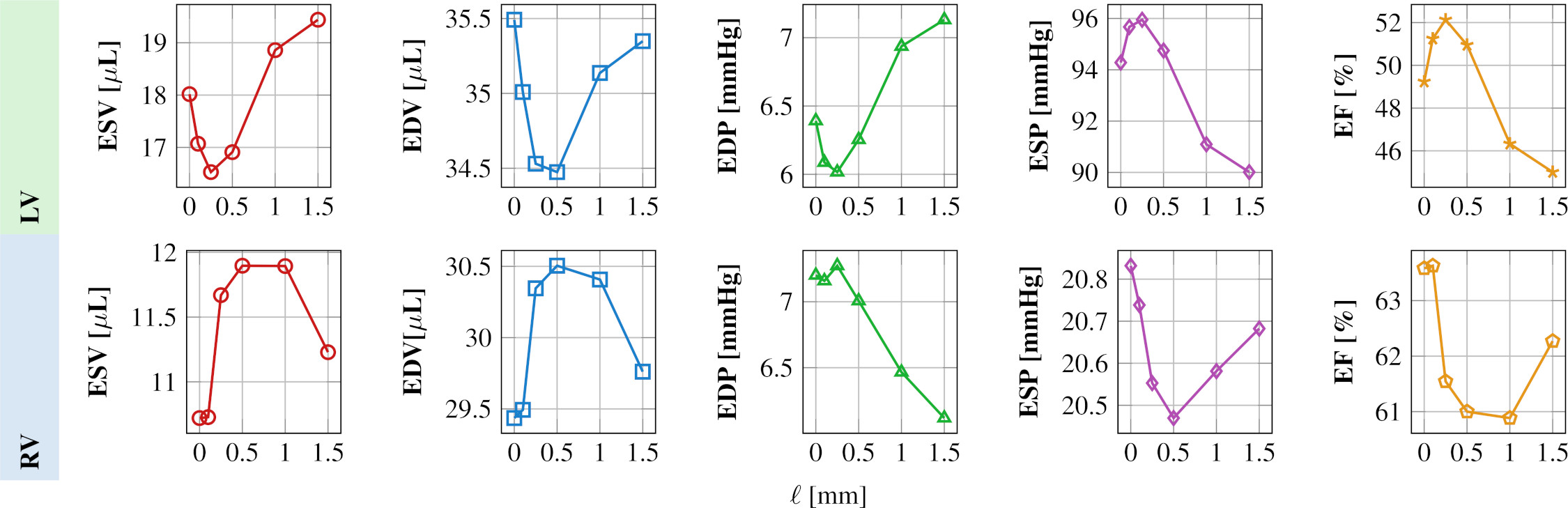}
    \caption{\textbf{\acrshort{qoi} at varying regularization radius.} Pressures and volumes of \acrshort{lv} and \acrshort{rv} during systole and diastole, and corresponding ejection fractions, obtained for different values of the regularization radius $\ell$.}
    \label{fig:QoIPlot}
\end{figure}

Fig~\ref{fig:QoIPlot} shows the variation of several \acrshort{qoi} with respect to the fiber regularization radius $\ell$: the End Diastolic Volume (EDV) and Pressure (EDP), the End Systolic Volume (ESV) and Pressure (ESP) and  
the Ejection Fraction (EF), which provides a measure of ventricular efficiency. Fig~\ref{fig:QoIPlot} shows a maximum EF of LV at $\ell = \SI{0.25}{\milli \meter}$ and reveals a non-monotonic effect of fiber regularization on cardiac contractility.
In addition, for \acrshort{lv}, the increases in EDP and ESP peaks at $\ell = \SI{0.25}{\milli \meter}$ indicate an increase in the amount of energy transferred to the blood.
This dual-phase behavior can be interpreted as the result of two distinct regimes induced by the fiber regularization radius $\ell$. 
For small values of $\ell$, increasing $\ell$ primarily reduces the local fiber disarray while preserving the underlying macroscopic fiber architecture. Moreover, \acrshort{lv} contraction becomes more effective, 
as myofibers act in a more coordinated manner, with less active mechanical energy dissipated along misaligned fibers.
Conversely, when the regularization radius $\ell$ approaches the characteristic length scale of the macroscopic fiber architecture, 
further smoothing progressively degrades the typical helical fiber organization, which is specifically structured to ensure a mechanically efficient contraction. 
As a result, the effectiveness of \acrshort{lv} contraction decreases. In addition EDV does not remain constant, when varying $\ell$, but it follows the same trend of the ESV, suggesting the effect of a residual tension that increases when $\ell \approx \SI{0.25}{\milli \meter}$. The residual tension limits the full expansion of the \acrshort{lv} at end diastole.

Based on the above observations, we interpret the fiber field obtained for $\ell = \SI{0.25}{\milli\meter}$ as a reliable representation 
of the effective macroscopic fiber architecture, while the discrepancy between this field and the originally measured one 
can be attributed to fiber disarray.
An opposite behavior with respect to \acrshort{lv} can be observed for \acrshort{rv}. The latter reflects the strong ventricular interdependence, with \acrshort{rv} dynamics predominantly dictated by LV activity \cite{milani2014small}.\\

\begin{figure}[thbp]
\includegraphics[width=\linewidth]{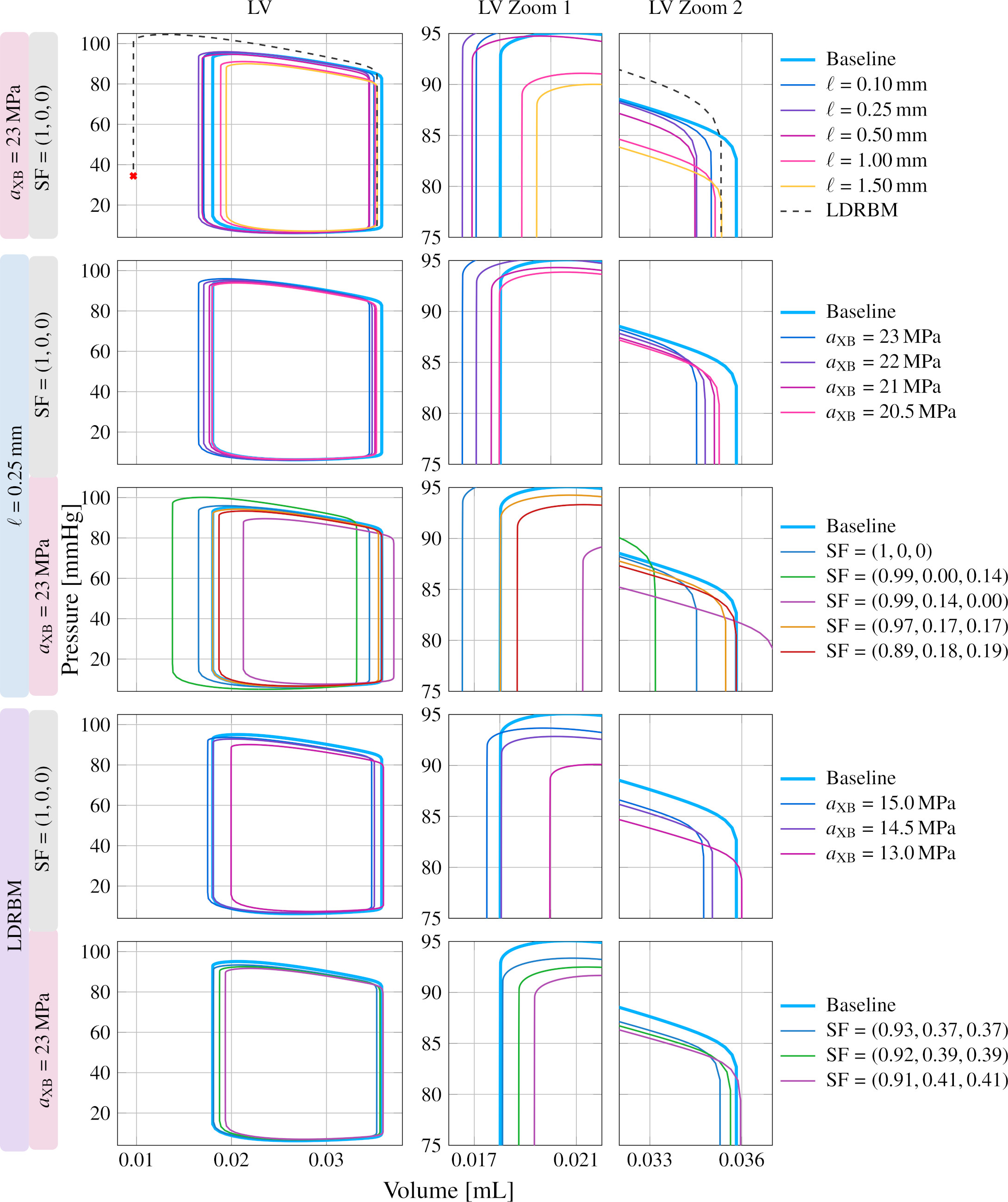}
\caption{\textbf{\acrshort{lv} PV loops comparison.}PV loops of \acrshort{lv} (first column) with zooms in end-systolic (\acrshort{lv} Zoom 1, second column) and end-diastolic (\acrshort{lv} Zoom 2, third column) zooms. The first row shows PV loops resulting from the baseline simulation (i.e., with the experimental fibers) compared with regularized and \acrshort{ldrbm} fibers. Rows two and three illustrate PV loops using the regularized fiber field, while rows four and five show PV loops with the \acrshort{ldrbm} fibers. The calibration was performed by adjusting $\CBstiffness$ (rows two and four) and $\SF$ (rows three and five).}
\label{fig:PVGrid}
\end{figure}

The first row of Fig~\ref{fig:PVGrid} compares the \acrshort{lv} PV-loop from the baseline simulation (using experimental fibers) with PV-loops obtained using regularized and \acrshort{ldrbm} fiber fields.
For $\ell \leq \SI{0.25}{\milli\meter}$, fiber regularization induces a leftward shift of the \acrshort{lv} PV loop. Further increasing the regularization radius $\ell$ causes the \acrshort{lv} PV loop to shift rightward, resulting in a decrease in \acrshort{ef} and, consequently, in pumping efficiency.
This behavior can be explained by the reduction of local fiber disarray when the regularization radius remains much smaller than the myocardial wall thickness, facilitating more a coordinated muscle contraction. Conversely, larger regularization radii $\ell$ progressively disrupt the original helical fiber architecture responsible for an optimal ventricular contraction. Fig~\ref{fig:PVGrid} also reports the PV loop obtained using \acrshort{ldrbm} fibers. In this case, the simulation failed during systole due to excessively high contraction levels, which caused the mesh element collapse. This behavior is attributed to the idealized nature of the rule-based fiber architecture, which enforces a perfectly helical fiber arrangement.

\subsubsection{Accounting for the presence of disarray in computational models}\label{subsubsec:Accounting}


Active force is largely directed along longitudinal directions, so fiber disarray can have a significant impact on both the magnitude and distribution of active tension.
We remark that most computational models account for this effect by introducing active stress components along directions orthogonal to the fibers \cite{eriksson2013influence, guan2020effect} (see Sec. 3). In Eq.~\eqref{eq:Piola}, this effect is represented by assigning nonzero values to $\mS$ and/or $\mN$. Commonly, however, these components are set to zero ($\mS = \mN = 0$), and the effect of disarray is implicitly included by calibrating the fiber contractility parameter \cite{berberouglu2014computational, rice2008approximate, land2012analysis}. 

In this study, the availability of a high-resolution experimentally measured fiber field, together with a modeling framework in which fiber disarray is progressively removed through smoothing, provides a unique opportunity to assess the impact of fiber disarray on active mechanics and to evaluate modeling strategies designed to account for it. Moreover, this framework enables us to assess whether the smoothed fiber field produces responses comparable to those observed with the experimental disordered architecture.

For the purposes of this analysis, and for the reasons outlined in Sec. 4.3.1, we use the smoothed fiber field obtained with $\ell = \SI{0.25}{\milli\meter}$ as representative of the macroscopic fiber architecture. We then surrogate fiber disarray using the proposed modeling strategies (i.e., by varying the \acrfull{sf} parameters, see Eq. 20) to assess their effectiveness in reproducing the functional consequences of fiber disarray.

We first investigate the impact of variations in the contractility parameter $\CBstiffness$. As shown in the second row of Fig~\ref{fig:PVGrid}, decreasing $\CBstiffness$ causes a rightward shift of the PV loops for the regularized fiber field ($\ell = \SI{0.25}{\milli \meter}$). Overall, adjusting contractility alone cannot fully reproduce the baseline LV response, but reducing $\CBstiffness$ from $\SI{23}{\mega\pascal}$ to $\CBstiffness = \SI{20}{\mega\pascal}$ (–13\%) provides a good approximation of the \acrfull{lv} behavior, see second row of Fig~\ref{fig:PVGrid}.


We then evaluate the effect of introducing active stress components along the cross-fiber directions $\myVec{s}_0$ and/or $\myVec{n}_0$. To isolate directional effects while preserving the overall level of activation,
we enforce the constraint $\mF^2 + \mS^2 + \mN^2 = 1$.
The third row of Fig~\ref{fig:PVGrid} illustrates that the sheet and normal components of active tension $\ActTension$ produce markedly different responses.
Adding normal-direction activation causes a leftward shift of the PV loop, while sheet-direction activation shifts it rightward. The latter effect is evident when comparing the fiber-normal activation case $\SF = (0.99, 0.00, 0.14)$ with the fiber-sheet activation $\SF = (0.99,\,0.14,\,0.00)$. When both the cross-fiber directions are activated, $\SF = (0.97,\,0.17,\,0.17)$, the resulting PV loop closely matches the baseline simulation (i.e., with the experimental fibers).

Finally, we consider a configuration in which active force is applied along the three directions $\myVec{f}_0$, $\myVec{s}_0$, and $\myVec{n}_0$, weighted according to the average projection of the experimentally measured fiber directions onto these orthogonal axes. Specifically, we compute the quantities
$\left| \myVec{f}^{\mathrm{meas}} \cdot \myVec{f}_0 \right|$, 
$\left| \myVec{f}^{\mathrm{meas}} \cdot \myVec{s}_0 \right|$, 
$\left| \myVec{f}^{\mathrm{meas}} \cdot \myVec{n}_0 \right|$,
where $\myVec{f}^{\mathrm{meas}}$ is the experimentally measured fiber direction, and then average them over the entire myocardium.
The resulting values correspond to the coefficients
$\SF = (0.89,\,0.18,\,0.19)$.
The use of absolute values reflects the directional invariance of the applied active force.
It is worth noting that $\sqrt{\mF^2 + \mS^2 + \mN^2} = 0.93 < 1$, so that this approach effectively
combines a reduction of the effective contractility with the introduction of active force
components in the cross-fiber directions. The resulting PV loop closely approximates the baseline simulation, exhibiting a better fit during
diastole compared to the case $\SF = (0.92,\,0.39,\,0.39)$, albeit with reduced agreement during systole.

The approaches previously applied to the experimental fiber field ( with
$\ell = \SI{0.25}{\milli\meter}$) are next tested using the \acrshort{ldrbm} fibers.
In this case, good approximations for the \acrshort{lv} PV-loop are obtained by reducing the crossbridge stiffness $\CBstiffness$, as shown in the fourth row of Fig~\ref{fig:PVGrid}, and also by 
introducing the misalignment for $\ActTension$, as illustrated in fifth row of Fig~\ref{fig:PVGrid}. The closest agreement is obtained for $\CBstiffness = \SI{12.5}{\mega \pascal}$ and for $\SF = (0.92, 0.39, 0.39)$. Nonetheless, the recovery of the baseline PV loop exhibits a larger approximation error relative to the regularized fiber field. This suggests that the average fiber orientation has a greater impact on contraction than fiber disarray.  As a result, replacing disarray with \acrshort{sf} does not fully capture the baseline contraction.

\subsubsection{Strain field analysis}\label{subsubsec:st_analysis}
To analyze the spatial distribution of deformation, we compute the first three invariants of the Green--Lagrange strain tensor $E(x)$
\begin{equation}
\myTens{E}(\mathrm{x}) = \frac{1}{2}(\myTens{F}^T(\myVec{x}) \myTens{F}(\myVec{x}) -\myTens{I}),
\end{equation}
where $\myTens{F}(\myVec{x})= I + \nabla \myVec{d}(\myVec{x})$ is the deformation gradient and $\myTens{I}$ is the Identity tensor. 
The first invariant $I_1$, given by the trace of the tensor $\myTens{E}$, $\myTens{I}_1=\operatorname{tr}(\myTens{E})$, provides a measure of the overall level of strain and is commonly associated with volumetric dilatation in the small-to-moderate deformation regime.
The second invariant $I_2=\frac{1}{2}[(\operatorname{tr}(\myTens{E})^2)-\operatorname{tr}(\myTens{E}^2)]$ captures distortional deformation by describing the deviatoric part of the strain tensor, thereby quantifying shape changes independently of volume variations.
The third invariant $I_3$, given by the determinant of the Green--Lagrange strain tensor $I_3=\det(\myTens{F}^T \myTens{F})$ characterizes higher-order nonlinear strain effects, and is associated with changes in volume.

\begin{figure}[H]
\centering
\includegraphics[width=0.75\linewidth]{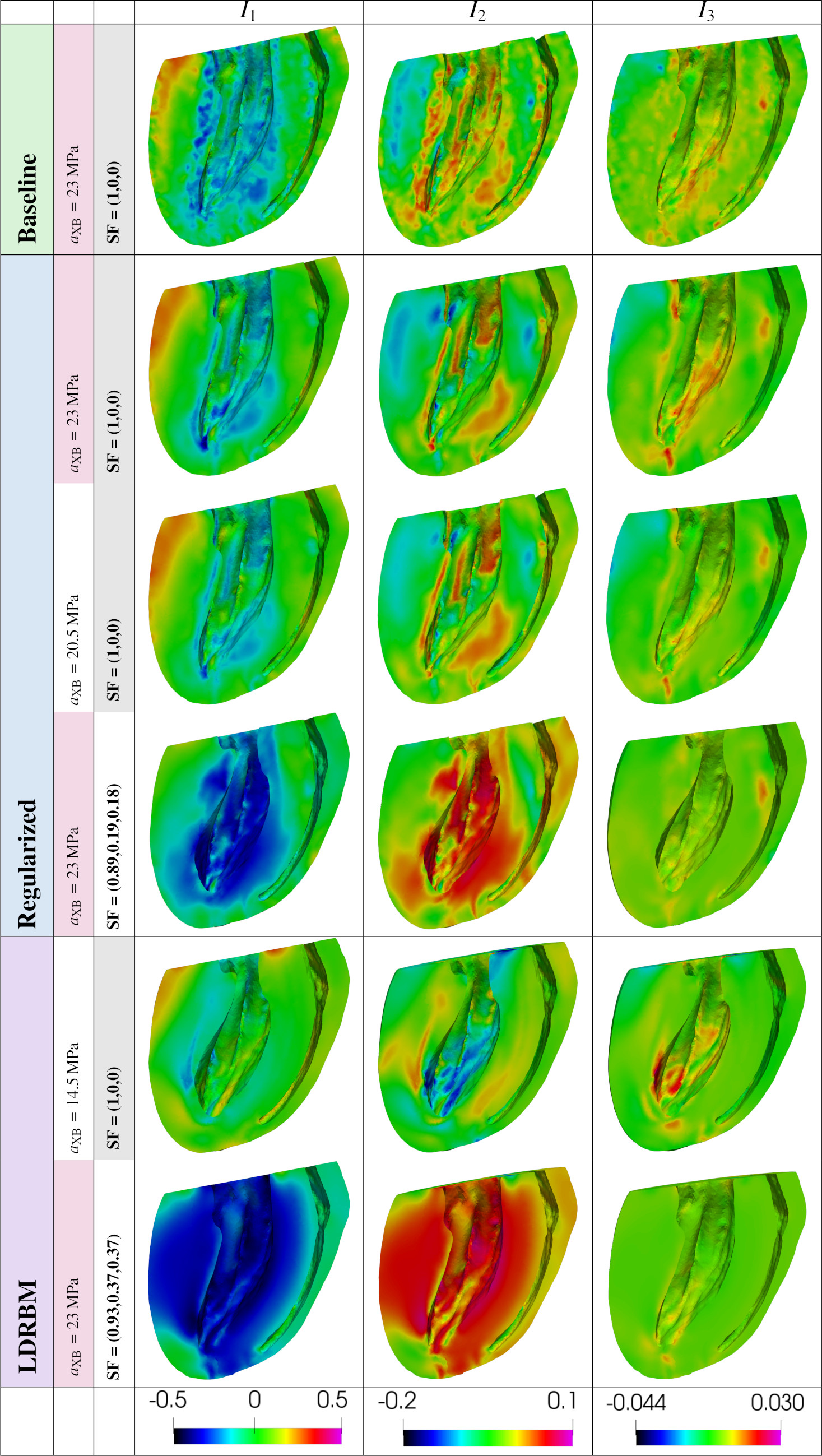}
\caption{\textbf{Invariants analysis.} Comparison of the first three invariants $I_1$, $I_2$, $I_3$ of the Green-Lagrange strain tensor $E$ (see Eq. X) obtained from different fiber fields: experimental (baseline), regularized fibers and \acrshort{ldrbm}.}
\label{fig:StrainFig}
\end{figure}

Fig~\ref{fig:StrainFig} compares the strain fields obtained for the baseline simulation, the regularized fiber field ($\ell = \SI{0.25}{\milli\meter}$), and the \acrshort{ldrbm}-based simulation. For each fiber architecture, we also include the configurations that best reproduce the PV loop of the baseline simulation.
For the regularized fiber field, these correspond to the cases with
$\CBstiffness = \SI{20.5}{\mega\pascal}$ and $\SF = (0.89,0.19,0.18)$, whereas for the rule-based fibers, the selected configurations use $\CBstiffness = \SI{12.5}{\mega\pascal}$ and $\SF = (0.92,0.39,0.39)$.

As shown in the second row of Fig~\ref{fig:StrainFig}, the regularized fiber field leads to a more diffused distribution of strain. While the main patterns of dilatation and distortion are preserved (with respect to the baseline), high-frequency spatial variations linked to fiber disarray are reduced.

A comparable strain pattern is observed in the reduced $\CBstiffness$ (see Fig~\ref{fig:StrainFig}, third row), with a reduction in the overall magnitude of both dilatation and distortional components. Conversely, the case with cross-fibers $\SF$ activation (see Fig~\ref{fig:StrainFig}, fourth row) exhibits marked discrepancies with respect to
the baseline simulation. 
Although strain concentration at the left endocardium is captured, local variations are largely not resolved.
In particular, dilatation at the left epicardium is absent, while spurious distortion appears in the
right myocardium.
This behavior can be attributed to the use of globally averaged $\SF$ parameters, which act
uniformly throughout the domain, whereas the experimentally fiber disarray shows significant heterogeneity.
As a result, the $\SF$ surrogate approach tends to enhance contraction uniformly, increasing contraction
at the endocardium while suppressing dilatation at the epicardium.

In the \acrshort{ldrbm} fiber simulation, with reduced $\CBstiffness$, the first invariant $I_1$ shows a contraction pattern similar to the experimental fiber cases, albeit with lower agreement than in the  regularized fiber field.
This suggests that the macroscopic fiber orientation is only globally captured by the rule-based approach. On the other hand, shear distortion is localized at the midwall, localized by the linear transmural variation of
fiber orientation, and at the left endocardium, where the sign of shear differs from the baseline.
Finally, the \acrshort{ldrbm} case with cross-fibers $\SF$ activation yields strain fields that are markedly different
from the baseline, with a predominantly negative first invariant and a uniformly positive second
invariant, thereby failing to reproduce the localized strain patterns observed in the baseline simulation.

\section{Discussion}\label{sec:Discussion}

In this work, we developed a biventricular electromechanical model of the murine heart, based on
high-resolution experimentally measured fiber orientations.
We introduced a strategy allowing to separate the macroscopic fiber architecture from local
fiber disarray, and applied it to investigate their respective roles in electromechanical
simulations.
The resulting responses were compared with those obtained using a surrogate fiber field generated with a \acrshort{ldrbm} approach \cite{doste2019rule}.
Finally, we assessed the influence of fiber architecture on the different components of the model,
including electrophysiology, passive mechanics, active contraction, and global hemodynamic indicators.

From the numerical results, we observed that the electrophysiology is only
marginally sensitive to the presence of fiber disarray, whose main effect is to promote a
more isotropic propagation of the activation wavefront (see Sec. 4.1). Moreover, the activation maps are well
approximated even with the rule-based fibers, indicating that \acrshort{rbms} capture the
macroscopic fiber architecture with sufficient accuracy to reproduce physiological activation
patterns. 
However, these considerations apply only to ventricular morphology under physiological conduction conditions, and can substantially change in the atria and/or under pathological conditions, where altered conduction properties or heterogeneous substrates arise due to scar tissue or fibrosis\cite{gander2023accuracy}.



The passive mechanical response is influenced by the macroscopic fiber architecture: replacing the experimental fiber field with the \acrshort{ldrbm} configuration results in a volume variation of approximately 4\% under purely passive loading at a pressure of $\SI{37.5}{\mega\pascal}$. Instead, the effect of microscopic fiber disarray -- represented in our analysis by increasing the smoothing parameter from $\ell = \SI{0}{\milli\meter}$ to $\ell = \SI{0.25}{\milli\meter}$ -- is much weaker in the passive regime, leading to volume changes of only about 0.4\%. For both macroscopic fiber architecture and microscopic disarray, the impact on passive mechanics is smaller than that observed in active contraction. A possible explanation is that active stress is primarily generated along the fiber direction and is therefore highly sensitive to fiber orientation. Conversely, passive elastic response is also supported by material stiffness in directions orthogonal to the fibers, resulting in a less pronounced anisotropy and, consequently, a reduced sensitivity to fiber orientation.

As a matter of fact, the greatest impact of fiber macroscopic architecture and microscopic disarray is observed in the electromechanical
simulations (see Sec. 4.3.1), where active tension plays a central role in myocardial contraction. The PV loops differ substantially among the baseline (i.e. with the experimental fibers), regularized, and \acrshort{ldrbm} cases.
In particular, a regularization radius around $\SI{0.25}{\milli\meter}$ yields the
maximum \acrshort{ef}, indicating that the reduction of microscopic disarray enhances the mechanical
efficiency of the pumping function. Further increasing the regularization radius, however, leads to
a reduction of the \acrshort{lv} \acrshort{ef}, suggesting that excessive alteration of the
macroscopic fiber organization impairs efficient ventricular contraction. Simulations performed with rule-based fibers and the same contractility parameter used for the experimental fibers
exhibit even stronger contraction, ultimately leading to numerical failure due to the highly
organized helical architecture.

We investigated two strategies commonly used in the cardiac modeling community to surrogate the effects of fiber disarray with an active stress approach (see Sec. 4.3.2): an effective reduction of contractility with respect to the pure longitudinal fiber contractile strength, and the introduction of active stress components in the
sheet and/or normal directions. In this analysis, simulations with experimental fibers were taken as the reference baseline, with the aim of validating these strategies as
practical modeling approaches to reproduce the behavior observed when high-resolution fiber data are available. The results showed that both strategies can yield pressure--volume loops close to the
baseline. In particular, the best match was obtained with a 13\% reduction in contractility, while when introducing active stress in orthogonal directions, the best agreement
is obtained by activating both the sheetlet and cross-fiber components with a ratio $\SF = (0.97,\,0.17,\,0.17)$. However, despite the good agreement at the organ level, significant discrepancies emerge at the tissue level. In particular, the inclusion of active stress in cross-fiber directions leads to strain fields that differ markedly from the baseline case. 
On the other hand, a simple effective reduction of
contractility, despite its phenomenological nature, is better able to reproduce the strain
distributions. These findings highlight an intrinsic limitation of modeling approaches based on the
addition of orthogonal active stress components. These results suggest that more faithful descriptions are likely to require
rigorous upscaling procedures that link microscopic fiber distributions to effective macroscopic
active stress formulations, which will be the subject of future work.

Finally, the regularized fiber fields consistently yield better agreement with the baseline results
compared to the \acrshort{ldrbm}, suggesting that an accurate representation of the macroscopic fiber
architecture plays a dominant role compared with the explicit modeling of microscopic disarray.

Previous studies have established the capability of \acrshort{rbms} to reproduce signal propagation in the ventricles 
\cite{bayer2005laplace,doste2019rule}. The analysis shown in the present work confirms these results and highlights the low impact of high frequency fibers misalignement in the signal propagation. Furthermore, our results suggest a negligible effect of fiber orientation on passive inflation. These findings are not comparable with those reported in \cite{palit2015computational}, where the authors performed a sensitivity analysis of fiber orientation angles on diastolic mechanics. They observed a clear dependence of myocardium elastic properties on fiber mean orientation, considering large variations of the fiber angles at endocardium and epicardium. Nevertheless, their results are derived from a \acrshort{rbm} and thus account solely for an idealized fiber architecture.

Regarding electromechanical simulations, our results highlight that a major impact of the  myofiber orientations and related disarray. 
The observed increase in \acrshort{edv} \cite{eriksson2013influence} for fibers oriented in the longitudinal direction, the enhancement of \acrshort{ef} in \acrshort{rbm} cases 
\cite{guan2020effect}, and the necessity to reproduce complete 3D fiber orientation for accurate strain distribution \cite{gil2019influence} are all consistent with literature findings \cite{gil2019influence,guan2020effect}. Additionally, the effects of active tension misalignment in sheet and normal directions are in agreement with previous results \cite{eriksson2013influence, guan2020effect}. In particular the active stress factor in the sheet direction decreases the \acrshort{ef}, while the active stress factor in the normal direction induces an opposite effect \cite{guan2020effect,piersanti20223d}.

Compared to previous studies, the present work represents, to the best of our knowledge, the first integration of fiber architecture data at this spatial resolution (\SI{96}{\micro \meter}) within a fully coupled electromechanical cardiac model. In addition, the proposed regularization framework enables a systematic decoupling of the effects of macroscopic fiber organization from the microscopic fiber disarray. This provides, for the first time, the opportunity to independently examine their functional roles and to critically evaluate commonly adopted modeling approaches for representing fiber disarray in electromechanical simulations.

A limitation of the present study is that the analysis was restricted to a single murine heart geometry. More definitive and quantitative conclusions will require extending the investigation to a larger cohort of geometries in order to account for inter-subject variability. Nevertheless, the electromechanical results presented here are robust and highlight clear trends, allowing us to identify important limitations of both rule-based fiber models and indirect modeling approaches for fiber disarray. Moreover, only physiological activation patterns were considered in this work. As a consequence, the conclusions drawn -- in particular regarding electrophysiology -- may not directly extend to pathological conditions, where altered conduction properties and heterogeneous substrates, such as scar or fibrosis, are expected to amplify the role of fiber architecture. Extending the present framework to pathological activation scenarios represents a relevant direction for future studies.

While direct quantitative translation to the human heart is not straightforward due to well-known inter-species differences in structure and function, the computational framework developed here is species-agnostic and may provide a basis for future investigations in larger animal models and, eventually, in human-specific geometries. From a clinical perspective, these findings would provide a coherent framework to interpret the role of myocardial fiber architecture in the cardiac function, and are also relevant for the construction of cardiac digital twins, where fiber architecture is typically approximated rather than directly measured. While capturing the macroscopic organization of myocardial fibers appears sufficient to reproduce global electrophysiological features, it can lead to significant inaccuracies in predicting mechanical function if microscopic disarray is not well captured, particularly when active contraction and strain distributions are of interest.

\section*{Authorship contributions}
Carlo Guastamacchia: Methodology, Software, Data Curation, Formal Analysis, Investigation, Visualization, Writing – Original Draft. Roberto Piersanti: Methodology, Software, Formal Analysis, Investigation, Writing – Review \& Editing. Francesco Giardini: Resources, Formal Analysis, Writing – Review \& Editing. Raffaele Coppini: Resources, Writing – Review \& Editing. Cecilia Ferrantini: Resources, Writing – Review \& Editing. Luca Dede': Writing – Review \& Editing. Leonardo Sacconi: Resources, Writing – Review \& Editing, Funding Acquisition. Francesco Regazzoni: Conceptualization, Methodology, Software, Formal Analysis, Investigation, Writing – Review \& Editing, Funding Acquisition

\section*{Acknowledgements}

This work is supported by ERC grant HeartCORE, funded by the European Union 770 (Grant Angreement 101198778); F.G. and L.S acknowledge Deutsche Forschungsgemeinschaft (DFG, German 782 Research Foundation) – Project 502822458.
C.G., L.D., and F.R. acknowledge the grant Dipartimento di Eccellenza 2023-2027, MUR, Italy.
C.G., R.P., L.D., and F.R. are members of GNCS, ``Gruppo Nazionale per il Calcolo Scientifico'' (National Group for Scientific Computing) of INdAM (Istituto Nazionale di Alta Matematica), Italy;
L.D. acknowledges the EuroHPC JU project dealii-X (grant number 101172493) funded under the HORIZON-EUROHPC-JU-2023-COE-03-01 initiative.

\bibliography{new_bibliography}

\newpage
\appendix
\section{Model calibration and baseline simulation setup}\label{app:Baseline_sim}

In this section, the baseline simulation and the corresponding parameters are detailed. Moreover, the calibration strategy for the electrophysiology, blood circulation, passive mechanics, and activation components is presented. First, we describe the employed experimental data, and then we illustrate the algorithm that we followed to calibrate all the components of the model, based on subject-specific activation maps, experimentally measured calcium and force traces, and on murine PV-loops observed in the literature (see \cite{tabima2010measuring,townsend2016measuring}).

\subsection{Baseline simulation setup}
As explained in Sec. 3, the electromechanical model is composed of electrophysiology, activation, mechanics and circulation core models. Electrophysiology is modeled through the eikonal equation, whose solution is the activation time. The initial conditions of the eikonal problem are composed of 5 stimulation sites activated at different times to reproduce the Purkinje network function \cite{stella2022fast}. The calcium ion concentration curve is then translated in time at each point of the domain by a quantity equal to the activation time computed through the eikonal equation. Mechanical activation is described by the RDQ20-MF model, which allows to compute the active tension $T_a$ through a system of ODEs \cite{regazzoni2020biophysically}. Then, $T_a$ enters in the definition of the active component of the Piola-Kirchhoff stress tensor $\myTens{P}$ as shown in Eq. (20). In the baseline simulation, the active tension is assumed to act only in the longitudinal fiber direction, therefore the $\SF$ are set as $\mF = 1$, $\mS = 0$, $\mN = 0$. We remark that $T_a$ depends linearly on $\CBstiffness$, which, in the baseline simulation, is set as $\CBstiffness = \SI{23}{\mega \pascal}$. Finally, the 0D model of blood circulation is coupled with the mechanics model to provide coherent boundary conditions \cite{piersanti20223d}.

The full set of parameters is listed in Tab.~\ref{tab:TabBaseline}.

\begin{table}
\includegraphics[width=\textwidth]{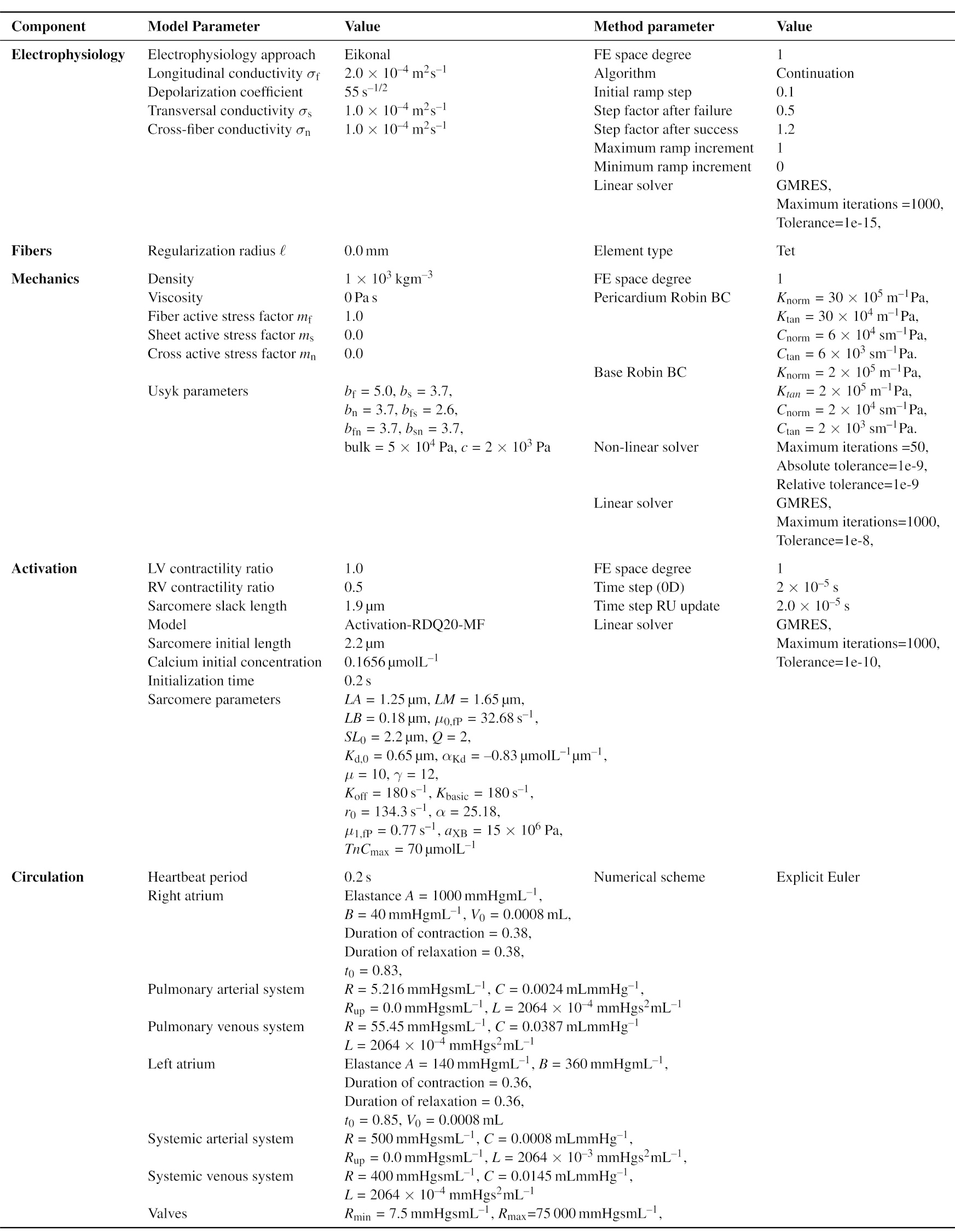}
\caption{Parameters for Baseline Model}
\label{tab:TabBaseline}
\end{table}

\subsection{Functional data for model calibration}

The functional data used for model calibration were acquired from one control  heart following the optical mapping protocol described in \cite{giardini2025correlative}. Briefly, the heart was excised from a heparinized and anesthetized male littermate control mouse, aged 28 weeks, cannulated at the aorta, and retrogradely perfused in a horizontal Langendorff configuration with oxygenated Krebs--Henseleit buffer at $\SI{36 \pm 1}{\degreeCelsius}$. Cardiac contraction was suppressed by adding blebbistatin to the perfusate. After stabilization of the electrocardiogram, the voltage-sensitive dye di-4-ANBDQPQ was injected into the aorta over $\SI{1}{\minute}$.

Transmembrane potential was imaged on the \acrshort{lv} free wall using wide-field fluorescence microscopy. Excitation was provided by a $\SI{625}{\nano \meter}$ LED filtered at $640/40~\si{\nano\meter}$, and the emitted fluorescence was collected through a $775/140~\si{\nano \meter}$ band-pass filter via a $\SI{685}{\nano \meter}$ dichroic beam splitter. Images were acquired on the central $128 \times 128$ pixels of an sCMOS camera (ORCA-Flash 4.0, Hamamatsu) at a frame rate of $\SI{1}{\kilo \hertz}$, yielding a pixel-size of $\SI{78}{\micro \meter} \times \SI{78}{\micro \meter}$ in the sample. The heart was electrically paced at the apex using an electrode connected to an isolated constant-voltage stimulator, delivering 15 pulses at $\SI{5}{\hertz}$.

Activation maps were derived from the recorded fluorescence signals. The mean baseline was subtracted and the signal was normalized on a pixel-by-pixel basis to yield the fractional fluorescence change ($\Delta F / F_0$). After signal inversion, a spatial binning of $4 \times 4$ pixels was applied, resulting in parameter maps with a spatial resolution of $\SI{312}{\micro \meter}$. Activation times were estimated using a cross-correlation-based approach: a reference pixel was selected, and the temporal shift of each pixel's fluorescence trace was computed relative to this reference, yielding a map of local activation delays. 

The two-dimensional activation map was co-registered with the three-dimensional mesoSPIM-based anatomical reconstruction of the same heart, following the morpho-functional registration pipeline described in \cite{giardini2025correlative}. The registered activation map, projected onto the corresponding anatomical mesh, served as the reference for calibrating the eikonal model parameters.

\subsection{Calibration procedure}
\begin{figure}
    \centering
    \includegraphics[width=\linewidth]{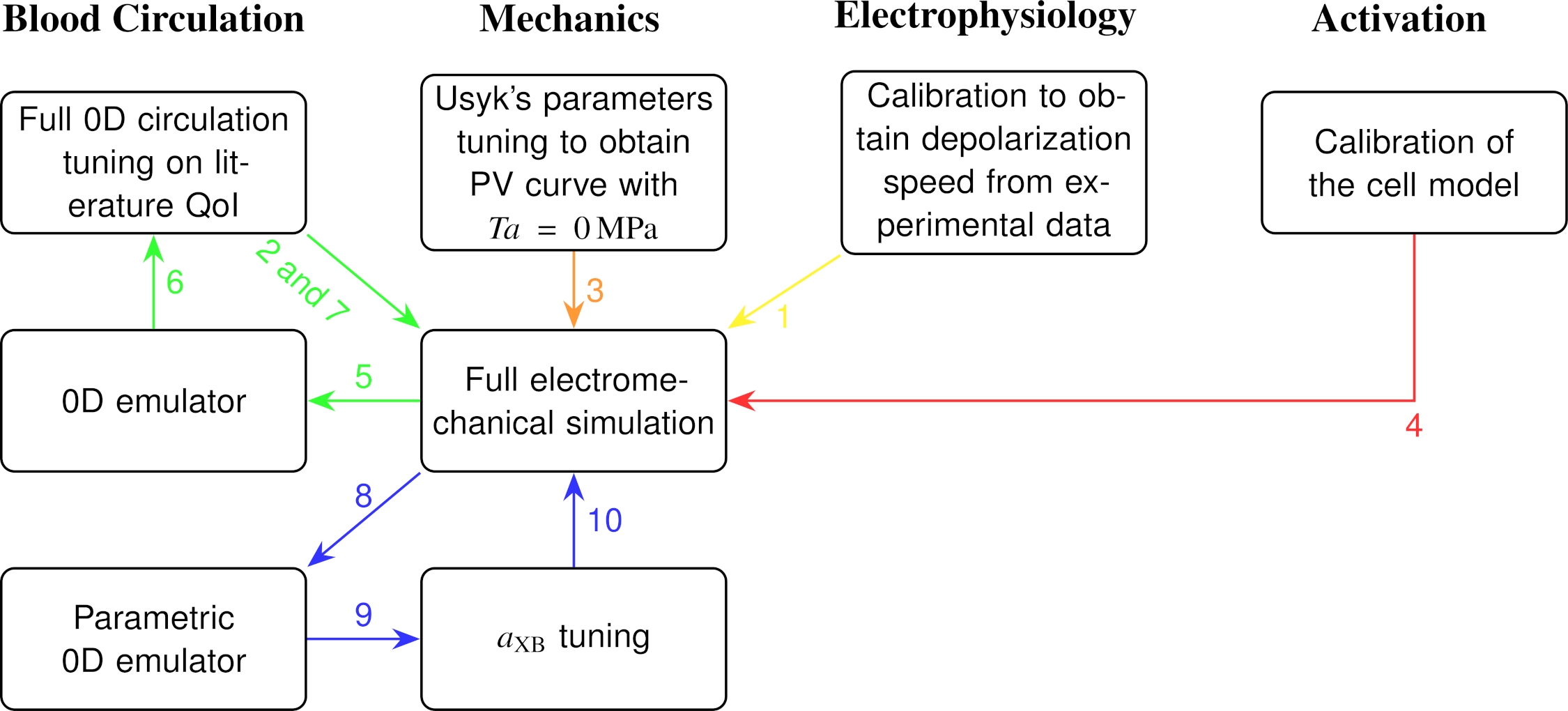}
    \caption{\textbf{Calibration workflow for the biventricular electromechanical model of the murine heart.} Arrows 1, 2, 3 and 4 represent the contribution of the four components in the electromechanical model. The green arrows (4, 5, and 6) depict a calibration loop that requires multiple iterations to calibrate the 0D circulation parameters. Arrows 7, 8, and 9 represent the final calibration loop for the crossbridge stiffness parameter.}
    \label{fig:TunMap}
\end{figure}

Fig.~\ref{fig:TunMap} shows a  representation of the calibration procedure within an iterative framework. The calibration procedure comprises multiple steps. The first step focuses on calibrating the electrophysiology parameters. The second step involves computing an initial guess for the set of parameters of the 0D circulation model, by temporarily surrogating the ventricles by means of time-varying elastance models. The third step consists of calibrating Usyk's model parameters in the passive mechanics model. Then the activation model is calibrated and a first guess for the crossbridge stiffness is computed. Afterward, the four core models are coupled in a comprehensive electromechanical simulation. From this model, a 0D emulator \cite{regazzoni2021accelerating} is built and the 0D circulation model calibration is fine-tuned together with the value of $\CBstiffness$. 

\subsubsection{Electrophysiology}
The eikonal model used in this work is calibrated based on an experimental activation map obtained through stimulation of the ventricular apex and imaging of the resulting fluorescence. From this measurement, the average depolarization speed in the mouse tissue is computed. Then, the conductivity parameters of the eikonal problem are adjusted  to recover the experimental depolarization speed (reported in Tab~\ref{tab:TabBaseline}).

\subsubsection{Blood-circulation}
The blood circulation model is calibrated based on ventricular PV-loops observed in literature \cite{tabima2010measuring,townsend2016measuring}. In the first step of the circulation calibration, the ventricles are modeled as elastances $E(t)$ varying in time t. Therefore, the pressure curve in each chamber can be computed as
\begin{equation}
    P(V,t) = E(t)(V - V_0),
    \label{eq:elastance}
\end{equation}
where $P$ is the pressure, $V$ is the blood volume and $V_0$ is the resting volume of the chamber. We solve multiple times the circulation model, coupled with Eq.~\eqref{eq:elastance}, to compute an initial guess of the circulation parameters. 

Following the calibration of the passive mechanics model, activation and electrophysiology, we use the 0D emulator proposed in \cite{regazzoni2021accelerating} to further calibrate the circulation parameters and speed up the convergence of the 3D model.
Specifically, the emulator can take advantage of 3D electromechanical simulations to build a 0D associated counterpart of the chamber (i.e. \acrshort{lv}) that reproduces the same PV-loop obtained in the 3D simulation. The latter is based on the assumption that the PV-loop of the chamber can be approximated by the weighted sum of the end-systolic pressure-volume relationship (ESPVR) and of the end-diastolic pressure-volume relationship (EDPVR). Therefore the ventricular pressure $P_{\mathrm{V}}$ can be written as
\begin{equation}
    P_{\mathrm{V}}^{\mathrm{emu}}(V,t) = (1 - \varphi_{\mathrm{act}}(t))P_{\mathrm{ED}}(V) + \varphi_{\mathrm{act}}(t)P_{\mathrm{ES}}(V)
    \label{eq: PV_loop_em}
\end{equation}
where $P_{\mathrm{ED}}(V)$ is the EDPVR, $P_{\mathrm{ES}}(V)$ is the ESPVR and the function $\varphi_{act}$ is a time-dependent function encoding the activation kinetics. Ideally, $\varphi_{act} = 0$, $\varphi_{\mathrm{act}} = 1$ corresponding to fully relaxed (at end-diastole) and fully contracted (at end systole) tissue respectively. $\varphi_{\mathrm{act}}$ is periodic in t with period $T_{\mathrm{HB}}$, that is the duration of a heartbeat. Note that inertia and damping effects are neglected within the 0D emulator. The procedure for building the emulator and reaching the limit cycle to be used in the 3D simulation consists of the following steps. 
\begin{enumerate}
    \item Run a full 3D simulation with an acceptable number of hearbeats (5 for this work).
    \item Build the 0D emulator of the 3D chamber using the PV-loop derived from the 3D simulation. 
    \item Run a 0D circulation simulation until the limit cycle is approached. During this step, the circulation parameters are calibrated in order to obtain a physiological PV-loop, based on those reported in the literature \cite{tabima2010measuring,townsend2016measuring}. This procedure is run in real time, avoiding running multiple simulations to calibrate the circulation parameters.
    \item Use the obtained PV-loops as initial conditions of the new 3D simulation.
\end{enumerate}
This method, which is called 3D-0D-3D V-cycle \cite{piersanti20223d}, is used to obtain the limit cycle of the calibrated circulation model in far less time than it would be required through multiple trials on the 3D simulations. The emulator, therefore, serves a dual role. Firstly, it speeds up the convergence process by skipping the computationally expensive transient behavior, enabling the 3D model to reach a steady state with far fewer simulated heartbeats. Moreover, it allows for efficient many-query scenarios, such as sensitivity analysis, parameter estimation, and uncertainty quantification, where the cardiac model is, in principle, multiple times for different parameter sets.

\subsubsection{Activation}
\begin{figure}[t]
    \centering
    \includegraphics[width=0.4\textwidth]{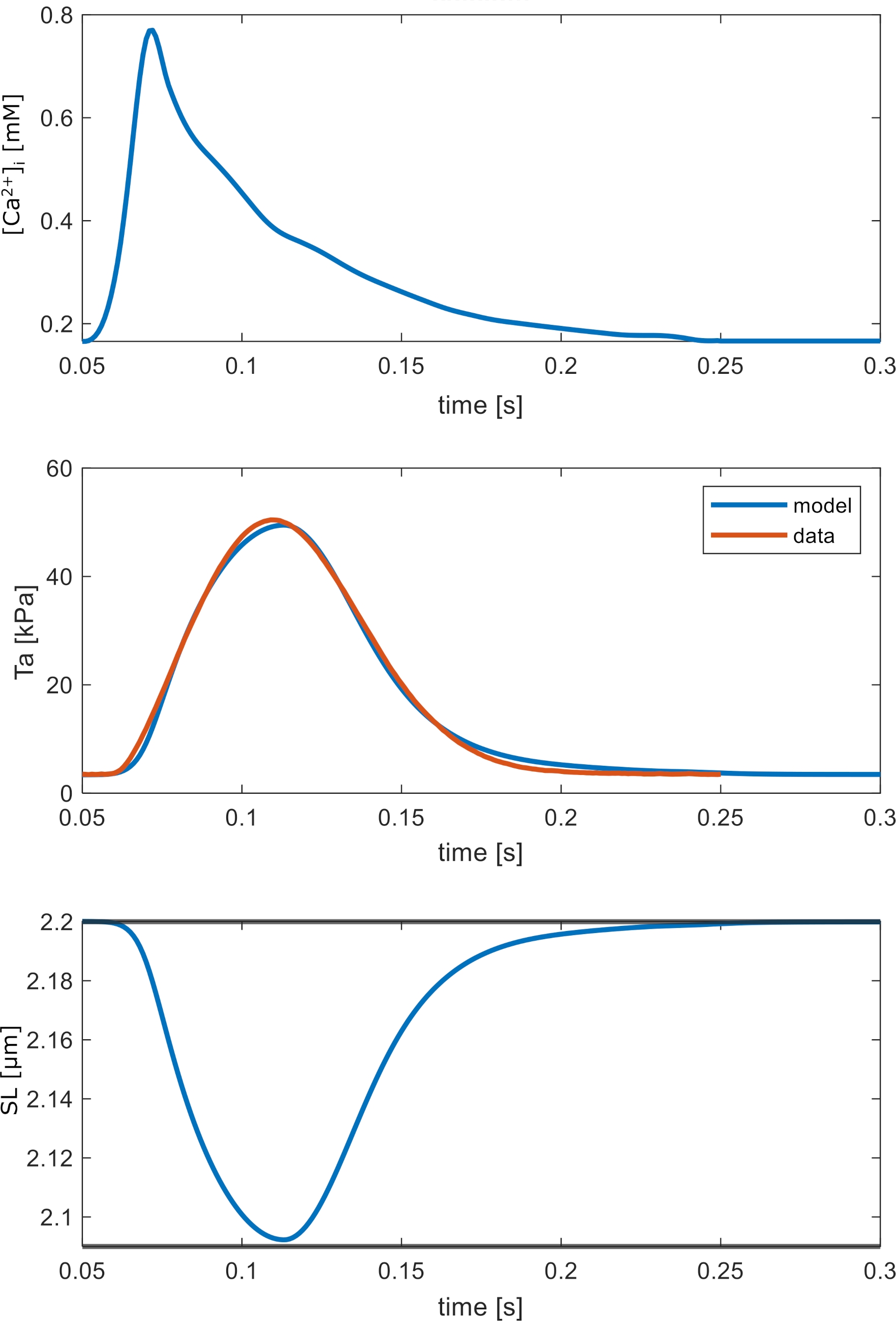}
    \caption{Calibration of the active force generation model for mouse cardiomyocytes. 
    Top: experimentally measured calcium transient, averaged over five consecutive cycles.
    Center: comparison between the experimental twitch force and the model prediction after calibration of the activation parameters. 
    The model reproduces the mechanical response by coupling the contractile element with a compliant elastic element adjusted to match the observed 10\% shortening at a diastolic sarcomere length of $\SI{2.2}{\micro\meter}$. 
    Bottom: sarcomere length transient predicted by the model under the same experimental conditions.}
    \label{fig:cell_model_calibration}
\end{figure}

The active force generation model RDQ20-MF was originally calibrated for rat and human cardiomyocytes \cite{regazzoni2020biophysically}. Therefore, its application to mouse cells requires a dedicated recalibration procedure. 
Since this process relies on detailed measurements obtained under specific conditions -- such as physiological temperature -- that were not available for the same animal model used here, we rely instead on data acquired from mice of comparable age at \SI{37}{\degreeCelsius}, with extracellular calcium concentration $\SI{1.8}{\micro M}$~\cite{ferrantini2016r4496c}.
Specifically, we use calcium transient measurements recorded during twitch contractions at \SI{5}{\hertz}. To reduce noise, the fluorescence signal is averaged over five consecutive cycles. 
To convert the Cal520 fluorescence signal $F(t)$ into intracellular calcium concentration, we use the standard relation $[\mathrm{Ca}^{2+}]_{\mathrm{i}}(t) = K_d \, \frac{F(t) - F_{\min}}{F_{\max} - F(t)}$ where $K_d = \SI{320}{\nano M}$ represents the dissociation constant of the dye, while $F_{\max} = 11110$ and $F_{\min} = 3600$ denote, respectively, the fluorescence levels at calcium saturation and at zero calcium.
The resulting calcium trace is then used to calibrate the RDQ20-MF model. To this end, we employed mechanical experiments performed on myofibers extracted from the same cell, consisting of shortening twitches characterized by a diastolic sarcomere length of $\SI{2.2}{\micro\meter}$ and a shortening amplitude of $10\%$. 
To model the experimental setup, we consider an elastic element in series with the contractile element, that takes into account the compliant portion of the sample. The stiffness of this elastic element is calibrated to match the observed shortening for a given developed force. 
Finally, the parameters of the activation model -- specifically those governing the kinetics ($k_{\mathrm{on}}$, $k_{\mathrm{basic}}$) and the cross-bridge contractility ($a_{\mathrm{XB}}$) -- are calibrated to minimize the mismatch between the predicted and measured force traces as shown in Fig.~\ref{fig:cell_model_calibration}.

\subsubsection{Mechanics}
For the passive mechanics, we start from Usyk’s parameter values obtained by the authors in \cite{nordbo2014computational}, where a murine \acrshort{lv} is inflated and the resulting pressure–volume curve is fitted by varying the stiffness parameters. Then, in order to match the literature passive pressure–volume relationship with the geometry used in the present study, which is characterized by a thin \acrshort{rv} wall, we vary the parameter $c$ in Usyk strain potential\cite{usyk2000effect}. 
Regarding the active part, to calibrate the crossbridge stiffness $a_{\mathrm{XB}}$, used in Eq. (21) a specific strategy is employed by using a combination of two emulators. This approach involves constructing two distinct emulators, each based on a different value for a particular parameter. Consequently, two separate full-scale simulations are used as datasets for emulator training. The two emulators are subsequently interpolated to create a parametric emulator, adding an extra degree of freedom. In particular given two different simulations, namely A and B, with $\CBstiffness = a_{\text{XB,A}}$ and $\CBstiffness = a_{\text{XB,B}}$ respectively, the parametrized pressure curve $p_{param}(V,t,\CBstiffness)$ becomes
\begin{equation}
    p_{\text{param}}(V,t,\CBstiffness) = \frac{\CBstiffness - a_{\text{XB,A}}}{a_{\text{XB,B}} - a_{\text{XB,A}}} p_{\mathrm{B}}(V,t) + \frac{\CBstiffness - a_{\text{XB,B}}}{a_{\text{XB,A}} - a_{\text{XB,B}}} p_{\mathrm{A}}(V,t).
    \label{eq:emu_parametric}
\end{equation}
We use this formulation to fit the PV-loops available in the literature \cite{tabima2010measuring,townsend2016measuring} by varying $\CBstiffness$.

Tab \ref{tab:TabBaseline} lists all the models and methods parameters used in the baseline simulation.
\section{Fiber distribution}\label{app:fibDist}
In this section, an analysis of fiber distribution in terms of mean angles and disarray is presented. 

\begin{figure}[t]
\includegraphics[width=\textwidth]{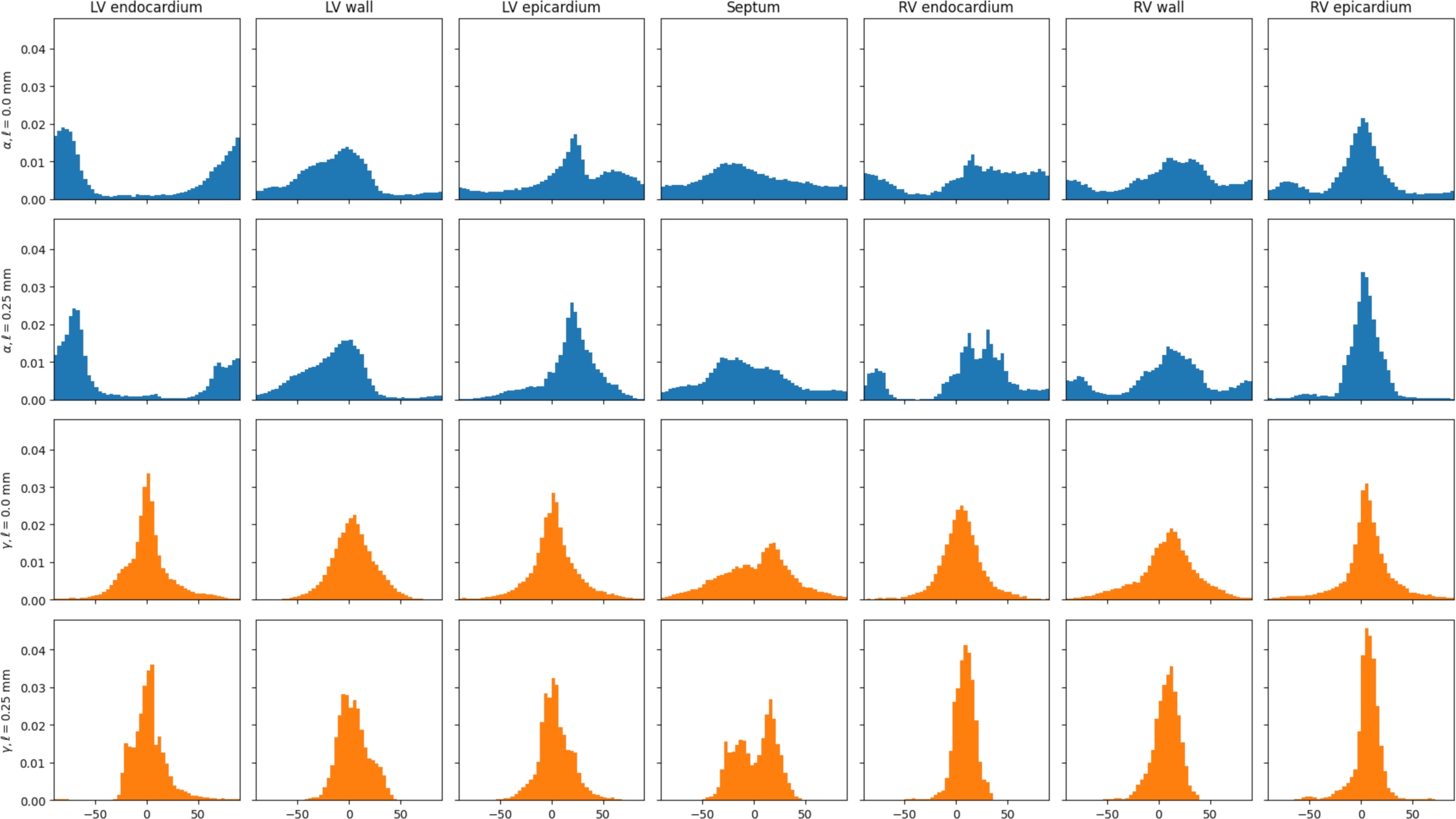}
\caption{\textbf{Distribution of the $\alpha$ and $\gamma$ angles in different zones of the domain.} Comparing the experimental fiber field ($\ell = \SI{0}{\milli \meter}$) with the regularized fiber field ($\ell = \SI{0.25}{\milli \meter}$).}
\label{fig:DistAlphaGamma}
\end{figure}

Fig.~\ref{fig:DistAlphaGamma} shows the distribution of the $\alpha$ and $\gamma$ angles in different regions of the biventricular geometry for the experimental fiber field ($\ell = \SI{0}{\milli \meter}$) and for the smoothed fibers ($\ell = \SI{0.25}{\milli \meter}$). The angle $\gamma$ exhibits mode values close to zero and its variance is much smaller with respect to $\alpha$. Furthermore, the smoothing process reduces the variance of the angles in the whole domain. In particular, the standard deviation is equal to $\SI{45.68}{\degree}$ for $\alpha$ and to $\SI{25.39}{\degree}$ for $\gamma$ in the experimental fiber field, while in the smoothed case the standard deviation decreases to $\SI{40.20}{\degree}$ for $\alpha$ and to $\SI{15.97}{\degree}$ for $\gamma$. The distribution of the angles varies largely across different regions of the myocardium.  The largest variation of the modal value of $\alpha$ is measured in the \acrshort{lv} passing from the endocardium (\SI{-76}{\degree}), through the wall (\SI{-1.8}{\degree}) to the epicardium (\SI{22}{\degree}), while in the zone of the \acrshort{rv} the variation is consierably lower (from \SI{5}{\degree} to \SI{14}{\degree}), given the thinner wall as shown in Fig.~\ref{fig:DistAlphaGamma}.

\begin{figure}[h]
\includegraphics[width=\textwidth]{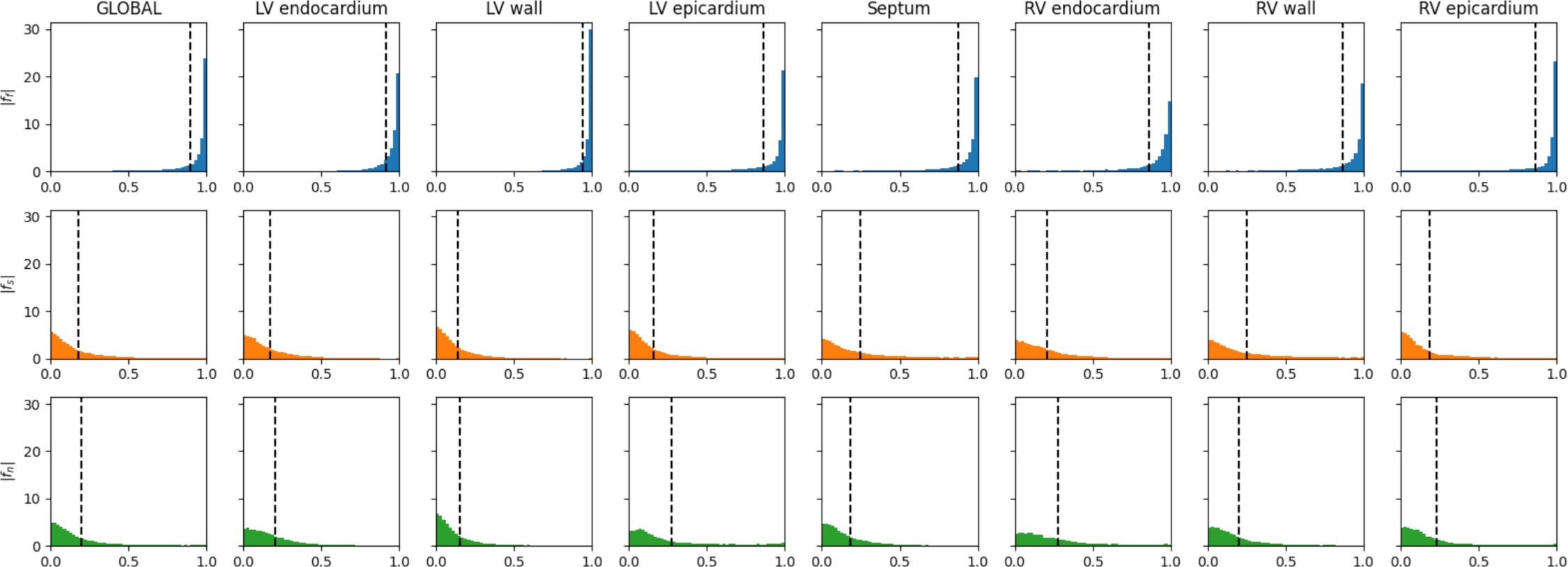}
\caption{\textbf{Distibution of $\hat{|\myVec{f}_f|}$, $\hat{|\myVec{f}_s|}$ and $\hat{|\myVec{f}_n|}$.}}
\label{fig:DistFibProjected}
\end{figure}

\begin{figure}[!h]
    \centering
    \includegraphics[width=\linewidth]{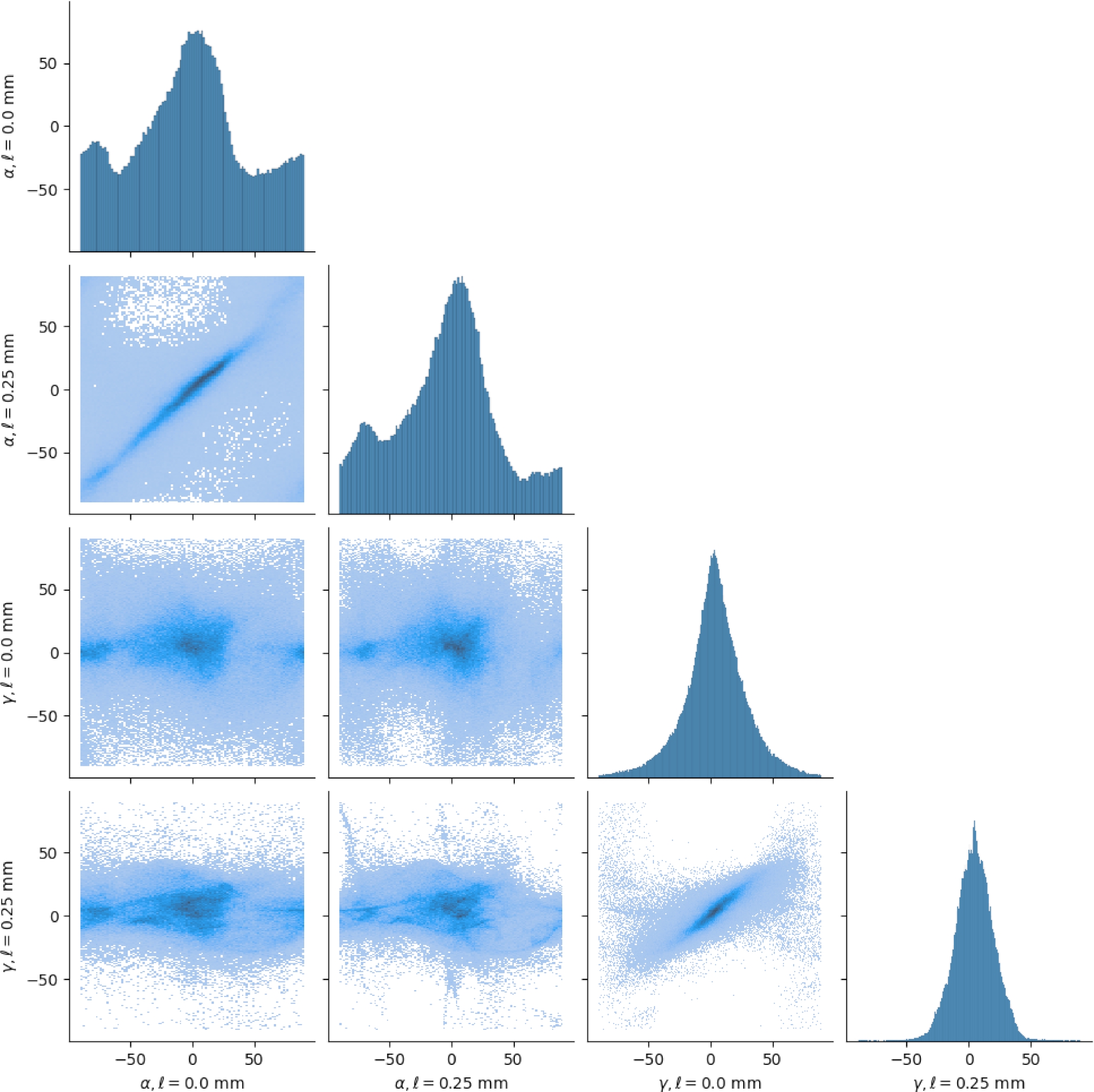}
    \caption{\textbf{Distribution of the angles $\alpha$ and $\gamma$ in the Experimental fiber field and in the smoothed fiber field.} Histograms showing angle distributions (diagonal boxes) and scatterplots showing the cross-correlations (off-diagonal boxes).}
    \label{fig:DistRel}
\end{figure}

 Fig.~\ref{fig:DistRel} shows the correlation of the $\alpha$ and $\gamma$ angles in the experimental field and in the smoothed one with $\ell = \SI{0.25}{\milli \meter} $. The plots show good correlation between the measured fibers and the smoothed one. In addition, we observe that the angles $\alpha$ and $\gamma$ are only weakly correlated between each other, as the angle $\gamma$ remains close to zero while $\alpha$ exhibits a much larger variation.

To assess the average fiber misalignment, the measured fiber field was projected onto the smoothed field $(\myVec{f}_0, \myVec{s}_0, \myVec{n}_0)$, resulting in three components: 
$\myVec{f}_f = \myVec{f}^{\mathrm{meas}} \cdot {\myVec{f}_0}$, 
$\myVec{f}_s = \myVec{f}^{\mathrm{meas}} \cdot {\myVec{s}_0}$ and 
$\myVec{f}_n = \myVec{f}^{\mathrm{meas}} \cdot {\myVec{n}_0}$. Fig, \ref{fig:DistFibProjected} shows the absolute value distribution of the three components in each zone of the domain. Given the means $\hat{|\myVec{f}_f|}$, $\hat{|\myVec{f}_s|}$ and $\hat{|\myVec{f}_n|}$, the distance $1 - \hat{|\myVec{f}_f|}$ and the values of $\hat{|\myVec{f}_s|}$ and $\hat{|\myVec{f}_n|}$ are measures of the myofiber disarray. In particular, $\hat{|\myVec{f}_f|}$ varies from a maximum of $0.92$ in the left endocardium to a minimum of $0.85$ in the right endocardium. In the whole myocardium we measure mean values $\hat{|\myVec{f}_f|}=0.89$,$\hat{|\myVec{f}_s|}=0.18$, $\hat{|\myVec{f}_s|}=0.19$. These coefficients are set in active stress factors in the configuration \acrshort{sf}$=(0.89, 0.19, 0.18)$ for the Piola-Kirchhoff stress tensor (see Eq. (20)).
\section{PV-loop analysis}\label{app:PVLoop}
\begin{figure}[!th]
    \centering
    \includegraphics[width=\linewidth]{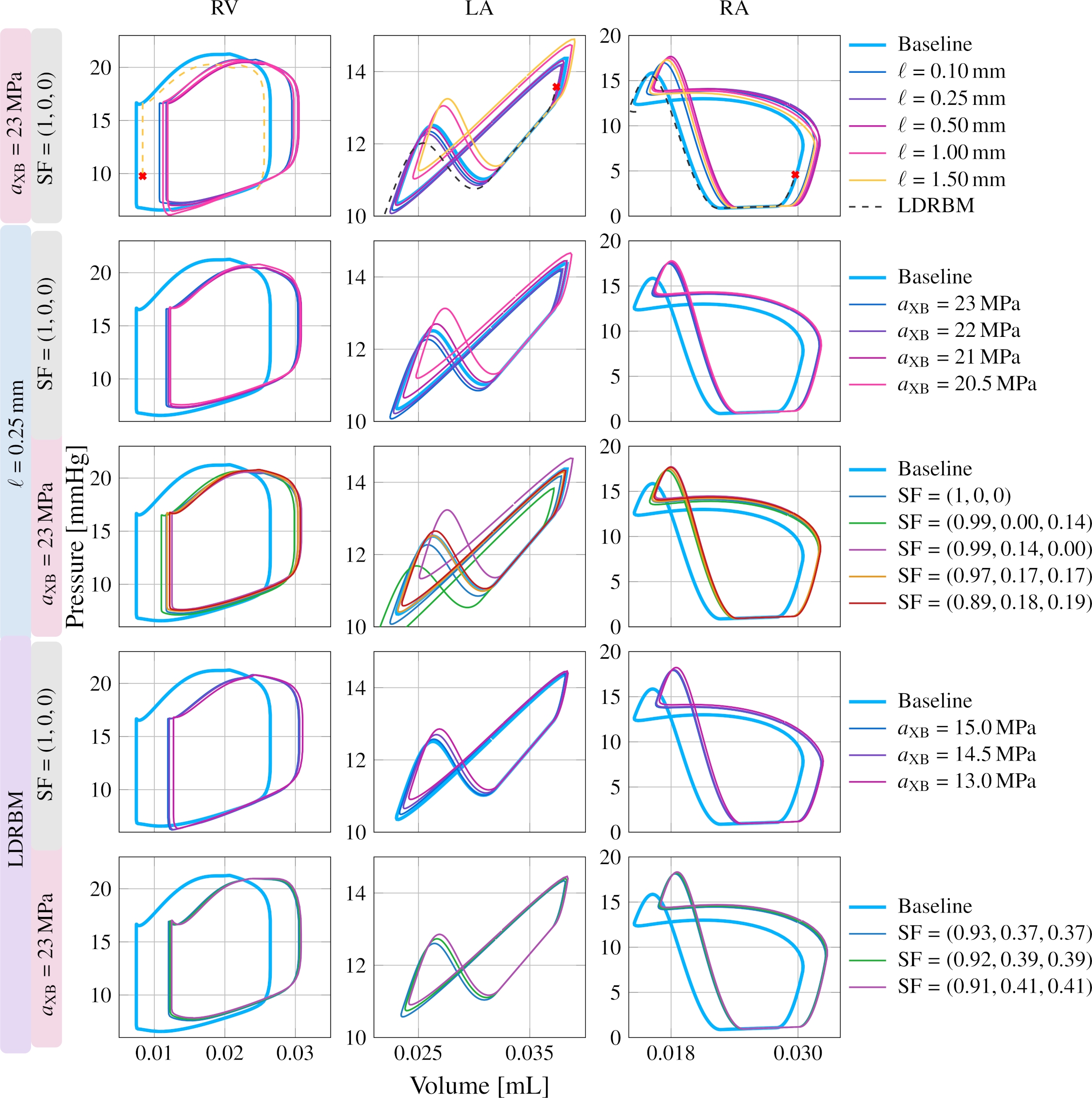}
    \caption{\textbf{PV-loops of (first column) \acrshort{rv}, (second column) \acrshort{la}, (third column) \acrshort{ra}.} (first Row) Comparison of PV-loops from experimental fiber field against regularized ones and the rule-based one, (second and third rows) recovery of the baseline PV-loop from regularized fiber field, (fourth and fifth rows) recovery of the baseline PV-loop from the rule-based one . The calibration is obtained varying (second and fourth rows) $\CBstiffness$  and (third and fifth rows) \acrshort{sf}.}
    \label{fig:GridAtrials}
\end{figure}

In this section we complement the analysis presented in Sec. 4.3.1 on \acrshort{lv} PV-loops with the PV-loops obtained by multiple simulations for \acrshort{rv}, \acrshort{la} and \acrshort{ra}. These plots are complementary to the \acrshort{lv} ones shown in Fig. 9. In particular, the first row of Fig.~\ref{fig:GridAtrials} shows the PV-loops obtained varying the regularization radius compared to the ones of the \acrshort{ldrbm}. In the second and third rows of Fig. \ref{fig:GridAtrials} we aim to recover the \acrshort{lv} PV-loops of the experimental fiber field and the regularized one with $\ell = \SI{0.25}{\milli \meter}$, varying the $\CBstiffness$ and \acrshort{sf}. Fig. \ref{fig:GridAtrials} in the fourth and fifth rows shows the PV-loops of the experimental fiber field and the \acrshort{ldrbm} ones. In the \acrshort{rv} we observe a shift to the right when increasing the regulrization radius up to $\ell = \SI{1.0}{\milli \meter}$, while the opposite trend is observed for $\ell=\SI{1.5}{\milli \meter}$. A similar behaviour is observed in \acrshort{ra} where the trend inversion is already visible for $\ell = \SI{1.5}{\milli \meter}$. We conclude that the smoothing of the fiber field does not causes an increase in the \acrshort{ra} and \acrshort{rv} efficiency. These results are consistent with the fact that the contraction of the \acrshort{rv} is strongly driven by \acrshort{lv}. Because the atria are represented by a 0D circulation model, the observed PV loops directly reflect the pressure and volume changes in the ventricles. The \acrshort{la} PV-loop follows the \acrshort{lv} behaviour with a left shift up to $\ell=\SI{0.25}{\milli \meter}$ and then a trend inversion. We recall that the \acrshort{ldrbm} simulation with \acrshort{sf} and $\CBstiffness$ equal to the baseline failed due to excessive contraction in the \acrshort{lv}. In the second and third rows we can observe that transversal and cross-fiber \acrshort{sf} activation and $\CBstiffness$ produce comparable effect on \acrshort{rv} and \acrshort{ra} with respect to \acrshort{lv}. In fact  the reduction of $\CBstiffness$ causes a shift to the right of the PV-loop for all the chambers, as the the activation of the sheet \acrshort{sf}, while the activation of cross-fiber \acrshort{sf} causes a shift to the left. However, the magnitude of this effect is considerably lower with respect to \acrshort{lv}. Passing from $\CBstiffness=\SI{23}{\mega \pascal}$ to  $\CBstiffness=\SI{20}{\mega \pascal}$ causes a reduction of \acrshort{ef} of \acrshort{rv} from 65\% 62\%, while passing from the \acrshort{sf}=(1,0,0) to the \acrshort{sf}=(0.99,0.14,0.00) reduces the \acrshort{ef} from 60\% to 58\%. Therefore while the \acrshort{lv} PV-loop get closer to the baseline, the \acrshort{rv} one shifts away from the baseline. Similar results are also observed in the fourth and fifth rows, where the \acrshort{rv} PV loops of the LDRBM are shifted to the right compared to the baseline simulation. This behaviour contrasts with that of the LV, where the LDRBM is more efficient than the baseline. These observations indicate that, with the tested methodologies, it is not possible to simultaneously recover the LV and RV PV loops using the LDRBM and the regularized fiber field. This confirms the need for improved methodologies that account for fiber disarray and overcome the lack of patient-specific data for fiber field reconstruction.

%
%

\end{document}